\documentclass[11pt]{article}
\usepackage{mathtools}

\usepackage{color}

\mathtoolsset{centercolon}
\usepackage{microtype}
\usepackage[usenames,dvipsnames,svgnames,table]{xcolor}
\usepackage{amsfonts,amsmath, amssymb,amsthm,amscd}
\usepackage[height=9in,width=6.5in]{geometry}
\usepackage{verbatim}

\usepackage[margin=10pt,font=small,labelfont=bf]{caption}
\usepackage{latexsym}
\usepackage[pdfauthor={ Tang},bookmarksnumbered,hyperfigures,colorlinks=true,citecolor=BrickRed,linkcolor=BrickRed,urlcolor=BrickRed,pdfstartview=FitH]{hyperref}
\usepackage{mathrsfs}
\usepackage{graphicx}
\usepackage{enumitem}

\usepackage[inline,nolabel]{showlabels}

\newtheorem{theorem}{Theorem}[section]

\newtheorem{lemma}[theorem]{Lemma}
\newtheorem{proposition}[theorem]{Proposition}

\newtheorem{conj}[theorem]{Conjecture}
\newtheorem{claim}[theorem]{Claim}
\newtheorem{example}[theorem]{Example}

\newtheorem{definition}[theorem]{Definition}
\newtheorem{remark}[theorem]{Remark}

\numberwithin{equation}{section}

\newcommand{\be}{\begin{equation}}
	\newcommand{\ee}{\end{equation}}
\newcommand\ba{\begin{align}}
	\newcommand\ea{\end{align}}

\newcommand{\notion}[1]{{\bf  \textit{#1}}}

\begin {document}
\author{
	Pengfei Tang\thanks{Department of Mathematical Sciences, Tel Aviv
		University.   
		Email: \protect\url{pengfeitang@mail.tau.ac.il}.}
}
\date{\today}
\title{Return probabilities  on nonunimodular transitive graphs}
\maketitle

\begin{abstract}
	Consider simple random walk $(X_n)_{n\geq0}$ on a transitive graph with spectral radius $\rho$. Let $u_n=\mathbb{P}[X_n=X_0]$ be the $n$-step return probability and $f_n$ be the first return probability at time $n$. It is a folklore conjecture that on transient, transitive graphs $u_n/\rho^n$ is at most of the order  $n^{-3/2}$. We prove this conjecture for  graphs with a closed, transitive, amenable and nonunimodular subgroup of automorphisms. 
 We also conjecture that for any transient, transitive graph  $f_n$ and $u_n$ are of the same order and the ratio $f_n/u_n$ even tends to an explicit constant. We  give some examples for which this conjecture holds.  For a graph $G$ with a closed, transitive, nonunimodular subgroup of automorphisms, we prove a weaker asymptotic behavior regarding to this conjecture, i.e., there is a positive constant $c$ such that $f_n\geq \frac{u_n}{cn^c}$. 
\end{abstract}

\section{Introduction and main results}\label{sec: intro}
\subsection{Local limit law of return probability}

Suppose $G=(V,E)$ is a  locally finite, connected, infinite graph with vertex set $V$ and edge set $E$. Let $(X_n)_{n\geq0}$  be a simple random walk on $G$ started from $x\in V$ and denote by $u_n(x):=\mathbb{P}_x[X_n=x]$  the $n$-step return probability.  In particular $u_0(x)=1$. The spectral radius $\rho$ of $G$ is $\rho:=\limsup_{n\to\infty}u_n(x)^{1/n}$, which doesn't depend on the choice of $x$ (for instance see \cite[Theorem 6.7]{LP2016}). Set $a_n(x):=\frac{u_n(x)}{\rho^n}$. When $G$ is (vertex)-transitive, the quantities $u_n(x)$ and $a_n(x)$ don't depend on $x$ and we simply write them as $u_n$ and $a_n$ respectively. 

A graph is called transient if a simple random walk on the graph is transient. It is  known that $\sum_{n=0}^{\infty}a_n<\infty$  for transient, transitive graphs; for instance see \cite[Theorem 7.8]{woess2000random}.  Since $a_{2n}$ is also decreasing in $n$ (using a Cauchy--Schwarz inequality as in the proof of Lemma 10.1 in \cite{woess2000random}), one has that $a_{2n}=o(\frac{1}{n})$. If some odd terms $a_{2k+1}>0$, then one still has that $a_n=o(\frac{1}{n})$ since $\frac{a_{2n+1}}{a_{2n}}\to 1$ (for instance see Lemma \ref{lem: ratio of u_n+1 over u_n tends to rho}). Hence for transient, transitive graph one always  has that $a_n=o(\frac{1}{n})$.  So what's more can one say about the asymptotic behavior of $a_n$ for such graphs? The following conjecture is folklore.
\begin{conj}\label{conj: a_n is at most n^(-3/2)}
	If a graph $G$ is  transient and transitive, then one has that
	\[
		a_n\preceq n^{-\frac{3}{2}}.
	\]
\end{conj}
Here for two functions $g,h:\mathbb{N}\to[0,\infty)$, we write $f(n)\preceq g(n)$ to denote that  there exists a constant $c>0$ such that $f(n)\leq cg(n)$ for all $n\geq0$. We write $f(n)\succeq g(n)$ if $g(n)\preceq f(n)$.  We write $f(n)\asymp g(n)$ if both $f(n)\preceq g(n)$ and $f(n)\succeq g(n)$ hold. We write $f(n)\sim g(n)$ if $\lim_{n\to\infty}\frac{f(n)}{g(n)}=1$.

  There are transient, non-transitive graphs such that $a_n(x)$ is bounded away from zero; for example see certain radial trees  in \cite{Gerl1984continued_fraction}. 

It is known that Conjecture \ref{conj: a_n is at most n^(-3/2)} holds for all transient, transitive, amenable graphs, for example $\mathbb{Z}^d\,(d\geq3)$. Let's briefly review this: for a transient, transitive and amenable graph $G$, its spectral radius $\rho$ equals $1$ (for example see Theorem 6.7 in \cite{LP2016}) and thus Conjecture \ref{conj: a_n is at most n^(-3/2)} becomes $u_n\preceq n^{-\frac{3}{2}}$ in this case. If $G$ has polynomial growth rate, i.e., $|B(x,r)|=O(r^\kappa)$ for some real number $\kappa>0$, then $G$ is roughly isometric to a Cayley graph (hence they have the same growth rate); see the discussion in the paragraph below the proof of Theorem 7.18 on page 265 of \cite{LP2016}. Furthermore there is an integer $d>0$ such that $|B(x,r)|=\Theta(r^d)$; for example see Theorem 5.11 of \cite{woess2000random}. Since $G$ is transient, one must have $d\geq 3$. Therefore $u_n\preceq n^{-\frac{d}{2}}\leq n^{-\frac{3}{2}}$ (for example see Corollary 14.5 of \cite{woess2000random}). If $G$ has superpolynomial growth rate, then for any $d>0$, $|B(x,r)|\succeq r^d$ by Theorem 5.11 of \cite{woess2000random}. Hence $u_n\preceq n^{-\frac{d}{2}}$ for all $d>0$, in particular  $u_n\preceq n^{-\frac{3}{2}}$.


It is also known that Conjecture \ref{conj: a_n is at most n^(-3/2)} hold for hyperbolic graphs \cite{Gouezel2014,GL2013RW_Fuchsian},  certain free products \cite{Candellero_etal2012RW_asymptotics_free_products,Cartwright1988RW_free_product,Cartwright_and_Soardi1986rw_free_products,Woess1986rw_free_produts} and certain Cartesian product \cite{Cartwright_and_Soardi1987local_limit_theorem_cartesian_products}. In particular for a regular tree $\mathbb{T}_{b+1}$ with degree $b+1\geq3$, it is known that $a_{2n} \sim \frac{1}{\sqrt{2\pi}}\cdot \frac{b+1}{2b}\cdot n^{-3/2}$.  However beyond the several cases just mentioned Conjecture \ref{conj: a_n is at most n^(-3/2)} is generally open for non-amenable Cayley graphs. 

Our first result is that Conjecture \ref{conj: a_n is at most n^(-3/2)}  holds for a certain family of transitive and  nonamenable graphs. See Section \ref{sec: examples} for specific examples which Theorem \ref{thm: main theorem on a_n for amenable, nonuni} applies to. 
\begin{theorem}\label{thm: main theorem on a_n for amenable, nonuni}
	If $G$ is a locally finite, connected graph with  a closed, transitive, amenable and nonunimodular subgroup of automorphisms, then $a_n\preceq n^{-\frac{3}{2}}$.
\end{theorem}

\subsection{First return probability}
Suppose $G=(V,E)$ is a locally finite, connected graph and $(X_i)_{i\geq0}$ is a simple random walk on $G$. For $x\in V$, the first return probability $f_n(x)$ is defined as $f_n(x):=\mathbb{P}_x[X_n=x,X_i\neq x,i=1,\cdots,n-1],n\geq 1$. We use the convention that $f_0(x)=0$. If $G$ is transitive, then  $f_n(x)$  doesn't depend on $x$ and we simply write it as $f_n$. 

Consider the generating functions $U(x,x|z)=\sum_{n=0}^{\infty}u_n(x)z^n$ and $F(x,x|z)=\sum_{n=0}^{\infty}f_n(x)z^n$. When the graph $G$ is transitive, we simply write $U(z)$ and $F(z)$ for $U(x,x|z)$ and $F(x,x|z)$. Since $\rho=\limsup_{n\to\infty}u_n(x)^{1/n}$, the radius  of convergence $r_U$ for $U(x,x|z)$ satisfies $r_U=\frac{1}{\rho}$. It is well known that $U(x,x|z)=\frac{1}{1-F(x,x|z)}$ for $|z|<\frac{1}{\rho}$; for instance see \cite[Lemma 1.13(a)]{woess2000random}.  Let $r_F$ be the radius of convergence of $F(x,x|z)$. It is known that $r_U=r_F$, in other words,
\begin{claim}\label{claim: same radii of conv for F(z) and U(z) for graphs}
		If $G=(V,E)$ is a locally finite, connected graph with spectral radius $\rho$, then for all $x\in V$,
		\[
		\limsup_{n\to\infty} f_n(x)^{1/n}=\rho.
		\]
\end{claim}
\begin{proof}
	This is a simple application of Pringsheim's theorem; for instance see Exercise 6.58 in \cite{LP2016}. 
\end{proof}

We conjecture something much stronger holds for all transient, transitive graphs.
\begin{conj}\label{conj: limit of the ratio of f_n over u_n}
	If $G$ is a locally finite, connected, transitive, transient graph, then 
	\[
		f_n\asymp u_n.
	\]

	Actually we conjecture the following equality holds:
	\[
	\lim_{n\to \infty,\mathsf{d}|n}\frac{f_n}{u_n}=\big(1-F(\rho^{-1})\big)^2\in(0,1),
	\]
	where $\mathsf{d}$ is the period of a simple random walk on $G$, i.e., $\mathsf{d}:=\gcd\{ n\geq 1\colon u_n> 0\}\in\{1,2\}$. 
\end{conj}

Conjecture \ref{conj: limit of the ratio of f_n over u_n} is known to hold for $\mathbb{Z}^d (d\geq3)$ (\cite{Doney_and_Korshunov2011}) and hyperbolic graphs (\cite[Proposition 4.1 and Theorem 1.1]{Gouezel2014}).   See Section \ref{sec: discussion of Conj 1.6} for more examples and  discussions on this.

Interestingly different behaviors occur for recurrent graphs. On $\mathbb{Z}$, it is well-known that $f_{2n}\sim \frac{1}{2\sqrt{\pi}n^{3/2}}$ while $u_{2n}\sim \frac{1}{\sqrt{\pi n}}$. On $\mathbb{Z}^2$, it happens that $f_{2n}\sim \frac{\pi}{n\log^2 n}$ (\cite{Jain_and_Pruitt1972} or \cite[Lemma 3.1]{Hamana2006pinned_RW}) while $u_{2n}\sim \frac{1}{\pi n}$. See \cite{Nachmias_and_Gurel2013nonconcetration} for some other results on first return probability on recurrent graphs.

The following is a partial result for nonunimodular transitive graphs, or more generally, graphs with a closed, transitive, nonunimodular subgroup of automorphisms. 
\begin{theorem}\label{thm: f_n over u_n is at least polynomial decay}
		If $G$ is a locally finite, connected graph with a closed, transitive and nonunimodular subgroup of automorphisms, then there is a constant $c>0$ such that 
		\[
		 f_n\succeq \frac{u_n}{n^c}.
		\]
		
\end{theorem}

Theorem \ref{thm: main theorem on a_n for amenable, nonuni} and Theorem \ref{thm: f_n over u_n is at least polynomial decay} are  also examples that sometimes nonunimodularity may help; see \cite[Theorem 1.2]{hutchcroft2017non} for another example on Bernoulli percolation.

\subsection{Organization of the paper and ideas of proof}
We prove Theorem \ref{thm: main theorem on a_n for amenable, nonuni} in Section \ref{sec: 2} and Theorem \ref{thm: f_n over u_n is at least polynomial decay} in Section \ref{sec: 3} respectively and then extend these results to the quasi-transitive case in Section \ref{sec: quasi-transitive case}. In Section \ref{sec: 5 examples for conjecture 1} we give some nonunimodular examples for Conjecture \ref{conj: a_n is at most n^(-3/2)}. Finally in Section \ref{sec: discussion of Conj 1.6} we discuss Conjecture \ref{conj: limit of the ratio of f_n over u_n} and give some examples for which this conjecture holds.

For Theorem \ref{thm: main theorem on a_n for amenable, nonuni} we observe that there is a natural choice of $\rho$-harmonic function $h$ and the Doob $h$-transform gives a new $p_h$-walk, and $a_n$ is the just the $n$-step return probability for this new $p_h$-walk. Next we observe that the $p_h$-walk is symmetric w.r.t.\ the level structure of the nonunimodular graph (Lemma \ref{lem: symmetric random walk on levels induced by the p_h walk}). Since one-dimensional symmetric random walk is well-understood, one can deduce that the probability that the  $p_h$-walk reaches a highest level $k$ and returns to level $0$ at time $n$ is  bounded by $(k\vee 1)^{3/2}/n^{3/2}$ (Lemma \ref{lem: probability of Y_n=0 and M_n=r}). Using the level structure again, on the event of reaching level $k$ and back to level $0$ at time $n$, the probability for the $p_h$-walk  returning to the starting point at time $n$ is bounded by $e^{-ck}$. Combining all this, we are done.

For Theorem \ref{thm: f_n over u_n is at least polynomial decay}, one can use mass-transport principle to deduce that the expected size of the intersection of a simple random walk path with level $k$ conditioned on returning at time $n$ is at most $ne^{-ck}$ (Proposition  \ref{prop: TMTP for the first moment intersect with a level}).  In particular, this implies that conditioned on returning at time $n$, the simple random walk has probability at least one half not reaching level $C\log n$ for large constant $C$. Then one can construct a first returning event as follows: first the walker starting from $x$ travels to a point $y$ in a lower level $-k=-C\log n$ with respect to $x$ in $k$ steps, and then does an excursion for  $n-2k$ steps without hitting the $k$-th level  with respect to $y$ (in particular not hitting $x$), and then travels back to $x$ in $k$ steps. This event itself has probability at least of order $\frac{u_n}{n^c}$ for some constant $c>0$.

\section{Proof of Theorem \ref{thm: main theorem on a_n for amenable, nonuni}}\label{sec: 2}

Suppose $G=(V,E)$ is a locally finite, connected graph. An automorphism of $G$ is a bijection $\phi:G\to G$ such that whenever $x$ and $e$ are incident in $G$, then so are the images $\phi(x)$ and $\phi(e)$. We denote by $\mathrm{Aut}(G)$ the group of automorphisms of $G$. Suppose $\Gamma\subset \mathrm{Aut}(G)$ is a closed subgroup of $G$, where we use the weak topology generated by the action of $\mathrm{Aut}(G)$ on $G$. We say $\Gamma$ is \textbf{transitive}, if for any pair of vertices $x,y\in V$, there is an element $\gamma\in\Gamma$ such that $\gamma(x)=y$. Denote by $x\sim y$ when $x,y$ are neighbors in $G$. (Recall that for two functions $f,g:\mathbb{N}\to(0,\infty)$, we also write $f(n)\sim g(n)$ if $\lim_{n\to\infty}\frac{f(n)}{g(n)}=1$. The meaning of the symbol $\sim$ can be easily determined from the context.)

\subsection{Amenability of graphs and groups}
Here we define of the amenability of graphs and groups. 
For a locally finite, connected graph $G=(V,E)$ and $K\subset V$, let $\partial_EK$ denote the edge boundary of $K$, namely, the set of edges connecting $K$ to its complement. 

\begin{definition}[Amenability of graphs]
	For a locally finite, connected, infinite graph $G=(V,E)$, let $\Phi_E$ be the edge-expansion constant given by 
	\[
	\Phi_E=\Phi_E(G):=\inf\Big\{ \frac{|\partial_EK|}{|K|};\emptyset\neq K\subset V\textnormal{ is finite} \Big \}.
	\] 
	We say the graph $G$ is amenable if $\Phi_E(G)=0$.
\end{definition}

A well-known result of Kesten states that for a locally finite, connected graph $G$, $G$ is amenable if and only if its spectral radius $\rho=1$; see \cite[Theorem 6.7]{LP2016} for a quantitative version.

\begin{definition}[Amenability of groups]
	Suppose $\Gamma$ is locally compact Hausdorff group and $L^{\infty}(\Gamma)$ be the
	Banach space of measurable essentially bounded real-valued functions on $\Gamma$ with respect
	to left Haar measure. We say that $\Gamma$ is amenable if there is an invariant mean
	on  $L^{\infty}(\Gamma)$.
	
	Here a linear functional on $L^{\infty}(\Gamma)$ is called a mean if it maps the constant
	function $\mathbf{1}$ to 1 and nonnegative functions to nonnegative numbers. Also a mean $\mu$ is called invariant if $\mu(L_\gamma f)=\mu(f)$ for all $f\in L^{\infty}(\Gamma),\gamma\in\Gamma$, where $L_\gamma f(x):=f(\gamma x),\forall\,x\in \Gamma$.
	
\end{definition}

The following theorem for transitive graphs is due to \cite{soardi1990amenability}; and Salvatori generalized it to the quasi-transitive cases.
\begin{theorem}[\cite{soardi1990amenability},\cite{Salvatori1992}]
	Let $G$ be a graph and
	$\Gamma$ be a closed quasi-transitive subgroup of $\mathrm{Aut}(G)$. Then $G$ is amenable iff $\Gamma$ is amenable and unimodular.
\end{theorem}

So in particular if a graph $G$ has a closed, transitive and nonunimodular subgroup of automorphisms, then $G$ is nonamenable.

For groups of automorphisms of graphs, Benjamini et al \cite{BLPS1999Group_inv} gave the following interpretation
\begin{lemma}[Lemma 3.3 of \cite{BLPS1999Group_inv}]\label{lem: amenability for automorphism group of graphs}
	Suppose $\Gamma$ is a closed subgroup of $\mathrm{Aut}(G)$ for the graph $G=(V,E)$. Then $\Gamma$ is amenable iff $G$ has a $\Gamma$-invariant mean.
	Here, a mean $\mu$ is  $\Gamma$-invariant  on $l^{\infty}(V)$ if $\mu(f)=\mu(L_\gamma f)$ for every $\gamma\in\Gamma,f\in l^{\infty}(V)$, where $L_\gamma f(x)=f(\gamma x)$ for $x\in V$.  
\end{lemma}

\subsection{Preliminaries on  nonunimodular transitive graphs}
Suppose $G$ is a graph and $\Gamma$ is a closed subgroup of $\mathrm{Aut}(G)$. 
There is a left Haar measure $|\cdot|$ on $\Gamma$ which is unique up to a multiplicative constant. We say $\Gamma$ is unimodular if the left Haar measure is also a right Haar measure; otherwise we say $\Gamma$ is nonunimodular. 

For a vertex $x$, denote by $\Gamma_x=\{\gamma\in \Gamma \colon \gamma(x)=x \}$  the stabilizer of $x$ in $\Gamma$. Let $m(x)=|\Gamma_x|$ be the left Haar measure of the stabilizer $\Gamma_x$.

\begin{lemma}[Lemma 1.29 of \cite{woess2000random}]\label{lem: uni of Haar msr via stabilizer}
	Suppose $\Gamma\subset \mathrm{Aut}(G)$ is a closed, transitive subgroup. For any $x,y\in V$, let  $\Gamma_x y$ denote the orbit of $y$ under $\Gamma_x$ and $|\Gamma_x y|$ denote the size of the orbit. Then 
\be\label{eq: ratio of haar msr of stabilizer in terms of ratio of the orbits}
\frac{m(y)}{m(x)}=\frac{|\Gamma_y x|}{|\Gamma_x y|}\,\,\,\forall\,x,y\in V.
\ee
\end{lemma}

\begin{proposition}[\cite{Trofimov1985}]\label{prop: unimodularity iff modular func is constant 1}
		Suppose $\Gamma\subset \mathrm{Aut}(G)$ is a closed, transitive subgroup. Then $\Gamma$ is unimodular if and only if 
		\[
		|\Gamma_y x|=|\Gamma_x y|\,\,\,\forall\,x,y\in V.
		\]
\end{proposition}

\begin{definition}\label{def: modular function}
	Suppose $\Gamma\subset \mathrm{Aut}(G)$ is a closed subgroup of automorphisms of the graph $G=(V,E)$.
	Define the modular function $\Delta:V\times V\to(0,\infty)$ by
	\[
	\Delta(x,y)=\frac{|\Gamma_y x|}{|\Gamma_x y|}.
	\]
\end{definition}

The following lemma contains the first two items in  \cite[Lemma 2.3]{hutchcroft2017non} that we shall need.
\begin{lemma}\label{lem: properties of modular func}
	The modular function $\Delta$ has the following properties.
	\begin{enumerate}
		\item[1.]  $\Delta$ is $\Gamma$-diagonally invariant, namely, 
		\[
		\Delta(x,y)=\Delta(\gamma x,\gamma y)\,\,\,\,\forall\,x,y\in V\,\,\forall\, \gamma\in\Gamma.
		\]
		
		\item[2.] $\Delta$ satisfies the cocycle identity, i.e.,
		\[
		\Delta(x,y)\Delta(y,z)=\Delta(x,z)\,\,\,\,\forall \, x,y,z\in V.
		\]
		
	\end{enumerate}
\end{lemma}

A key technique is the mass-transport principle. 
\begin{proposition}[Theorem 8.7 of \cite{LP2016}]\label{prop: TMTP}
		Suppose $\Gamma\subset \mathrm{Aut}(G)$ is a closed, transitive subgroup. 
	If $f:V\times V\to[0,\infty]$ is a $\Gamma$-diagonally invariant function, then 
	\be\label{eq: TMTP form 1}
	\sum_{v\in V}f(x,v)=\sum_{v\in V}f(v,x)\Delta(x,v)
	\ee
\end{proposition}

The following lemma is a simple application of the mass-transport principle. 
\begin{lemma}\label{lem: ratio of m(x) for neighboring vertices}
	Suppose $\Gamma$ is a closed, nonunimodular, transitive subgroup of $\mathrm{Aut}(G)$. 
	Write $B:=\big\{\frac{m(y)}{m(x)}\colon y\sim x \big\}$ for the set of all possible values of the modular function on two neighboring vertices. Write  $B_+:=\{q\in B\colon q>1\}=\{q_1,\cdots,q_k\}$ and  $B_-:=\{ q\in B\colon q<1\}$.  For $q\in B$ write  $t_q:=|\{y\colon y\sim x,\frac{m(y)}{m(x)}=q \}|$ for the number of neighbors of $x$ such that the  modular function $\Delta(x,y)$ takes the value $q$. Then 
	\begin{enumerate}
		\item[(i)]   $B_-=\{q^{-1}\colon q\in B_+\}$ and 
		\item[(ii)]   $t_{q^{-1}}=qt_q$ for all $q\in B$.
	\end{enumerate}
\end{lemma}
\begin{proof}
	For $q\in B$, define $f:V\times V\to (0,\infty)$ by 
	\[
	f(x,y):=\mathbf{1}_{\{y\sim x, \frac{m(y)}{m(x)}=q \}}.
	\]
	Obviously $f$ is $\Gamma$-diagonally invariant. 
	
	By the  mass-transport principle (Proposition \ref{prop: TMTP}),
	\[
	\sum_{y\in V}f(x,y)=\sum_{y\in V}f(y,x)\Delta(x,y),
	\]
	i.e.,
	\be\label{eq: t_q and t_q^-1}
	t_q=t_{q^{-1}}q^{-1}.
	\ee
	In particular, one has $t_q>0$ iff $t_{q^{-1}}>0$. Hence $B_-=\{q^{-1}\colon  q\in B_+\}$. Moreover, \eqref{eq: t_q and t_q^-1} gives the conclusion (ii).
\end{proof}

\subsection{A $\rho$-harmonic function}

Suppose $G=(V,E)$ is a transitive, locally finite, infinite graph with spectral radius $\rho$. 
Let $P$ denote the transition operator associated with simple random walk $(X_i)_{i\geq0}$ on $G$, i.e.,
\[
(Pf)(x)=\mathbb{E}_x[f(X_1)]=\sum_{y\in V}p(x,y)f(y),
\]
where $p(x,y)=\mathbb{P}_x[X_1=y]$. We also write $p^{(n)}(x,y)=\mathbb{P}_x[X_n=y]$ for the $n$-step transition probability from $x$ to $y$. In particular $u_n=p^{(n)}(x,x),\forall x$.
We say a function $f:V\to\mathbb{R}$ is \textbf{$\rho$-harmonic} if $Pf=\rho f$.

If there is a $\rho$-harmonic positive function $h$ on $V$, then one can define the Doob transform $p_h:V\times V\to(0,\infty)$ by
\[
p_h(x,y)=\frac{p(x,y)h(y)}{\rho\cdot h(x)}.
\]
Since $h$ is $\rho$-harmonic, the function $p_h$ defines a transition probability on $G$ and we call the corresponding Markov chain the $p_h$-walk. Recall that $a_n:=\frac{u_n}{\rho^n}$.
For any vertex $x$ of the transitive graph $G$, obviously the $n$-step transition probability of the $p_h$-walk satisfies:
\be\label{eq: n-step transition prob of p_h is just a_n}
p_h^{(n)}(x,x)=\frac{p^{(n)}(x,x)h(x)}{\rho^n h(x)}=\frac{u_n}{\rho^n}=a_n. 
\ee

\begin{lemma}\label{lem: rho-harmonic func for amenable nonuni}
		Let $G$ be a connected graph with  a closed, transitive, amenable and nonunimodular subgroup $\Gamma$ of automorphisms. Let $h:V\to(0,\infty)$ be given by $h(x)=\sqrt{m(x)}$.
	 Then the function
	$h$ is $\rho$-harmonic on $G$ and $\Gamma$-ratio invariant in the sense that 
	\[
	\frac{h(\gamma y)}{h(\gamma x)}=\frac{h(y)}{h(x)}\,\,\,\,\forall\,x,y\in V, \forall\, \gamma\in \Gamma.
	\]
\end{lemma}

For this lemma we need Theorem 1(b) from \cite{soardi1990amenability}. It says that if $G$ is a connected, transitive graph with spectral radius $\rho$ and degree $d$, and $\Gamma$ is a closed subgroup of $\mathrm{Aut}(G)$ which acts transitively on $G$, then one has that 
\[
\rho\leq \frac{1}{d}\sum_{y\colon y\sim x}\sqrt{\frac{|\Gamma_yx|}{|\Gamma_xy|}},
\]
with equality holds if and only if $\Gamma$ is amenable. 

\begin{proof}[Proof of Lemma \ref{lem: rho-harmonic func for amenable nonuni}]
 The $\Gamma$-ratio invariance of $h$  follows from the $\Gamma$-diagonally invariance of the modular function $\Delta$; see Lemma \ref{lem: properties of modular func}.
	
Since $\Gamma$ is amenable, Theorem 1(b) of \cite{soardi1990amenability} then implies the $\rho$-harmonicity of $h$: 
	\be\label{eq: rho-harmonic for square root of m(x) in case of amenable and nonuni}
	\rho=\frac{1}{d}\sum_{y\colon y\sim x}\sqrt{\frac{|\Gamma_yx|}{|\Gamma_xy|}}=\frac{1}{d}\sum_{y\colon y\sim x}\sqrt{\frac{m(y)}{m(x)}},
	\ee
	where $d$ is the degree of $G$ and the second equality is due to Lemma \ref{lem: uni of Haar msr via stabilizer}.
\end{proof}

\begin{proposition}\label{prop: properties of p_h}
	Given a positive, $\Gamma$-ratio invariant, $\rho$-harmonic function $h$ on a transient, transitive graph $G$, the $p_h$-walk  on $G$ is  transient, $\Gamma$-invariant and reversible. 
\end{proposition}
\begin{proof}
	Transience follows from \cite[Theorem 7.8]{woess2000random}: $\sum_{n\geq0}p_h^{(n)}(x,x)\stackrel{\eqref{eq: n-step transition prob of p_h is just a_n}}{=}\sum_{n\geq0} a_n<\infty $.

	Since $h$ is $\Gamma$-ratio invariant, $p_h$ is $\Gamma$-invariant:
	\[
	p_h(\gamma x,\gamma y)=p_h(x,y),\forall\,x,y\in X,\gamma\in\Gamma.
	\]
	
	Reversibility: let $\pi(x)=h(x)^2$, then 
	\be\label{eq: conductance for p_h}
	\pi(x)p_h(x,y)=h(x)^2\frac{p(x,y)h(y)}{\rho\cdot h(x)}=\frac{\mathbf{1}_{\{x\sim y \}} h(y)h(x)}{d\cdot \rho},
	\ee 
	where $d$ is the degree of $G$.
	Hence $\pi(x)p_h(x,y)=\pi(y)p_h(y,x)$. 
\end{proof}

\subsection{Proof of Theorem \ref{thm: main theorem on a_n for amenable, nonuni}}

Throughout this subsection, we assume $G$ is a connected graph with  a closed, transitive, amenable and nonunimodular subgroup $\Gamma$ of automorphisms. 

We first study the $p_h$-walk associated with the $\rho$-harmonic function $h(x)=\sqrt{m(x)}$ from Lemma \ref{lem: rho-harmonic func for amenable nonuni}. This random walk is a special case of the so-called ``square-root biased" random walk in \cite[Definition 5.6]{Tang2019heavy}.

Let $(S_n)_{n\geq0}$ be a $p_h$-walk on $G$ started with $o$. Let $(Y_n)_{n\geq0}$ be given by $Y_n:=\log \Delta(S_0,S_n)$. Then using the cocycle identity for the modular function (Lemma \ref{lem: properties of modular func}), we know that the increment sequence $(Z_i)_{i\geq1}$ is  a sequence of i.i.d.\ random variables, where $Z_i:=Y_i-Y_{i-1}=\log\Delta(S_{i-1},S_i),i\geq1$.
\begin{lemma}\label{lem: symmetric random walk on levels induced by the p_h walk}
The random walk	$(Y_n)_{n\geq0}$ is a symmetric random walk on $\mathbb{R}$ starting from $0$ with i.i.d.\ increments and the increments are bounded and have  mean $0$.
\end{lemma}
\begin{proof}
	From Lemma \ref{lem: ratio of m(x) for neighboring vertices}, the range of $Z_i$ is the finite set $\{ \log q\colon q\in B \}$. In particular, the increments are bounded.
	
	Notice that 
	\[
	\mathbb{P}[Z_1=\log q]=\sum_{y\colon y\sim x,\frac{m(y)}{m(x)}=q}p_h(x,y)=\sum_{y\colon y\sim x,\frac{m(y)}{m(x)}=q}\frac{1}{d}\cdot \frac{\sqrt{m(y)}}{\rho\sqrt{m(x)}}=\frac{t_q\cdot \sqrt{q}}{d\rho}. 
	\]
	In particular, $(Y_i)_{i\geq0}$ is symmetric: for any $q\in B_+$,
	\[
	\mathbb{P}[Z_1=\log q]=\frac{t_q\cdot \sqrt{q}}{d\rho}\stackrel{\eqref{eq: t_q and t_q^-1}}{=}\frac{t_{q^{-1}}\cdot \sqrt{q^{-1}}}{d\rho}=\mathbb{P}[Z_1=-\log q]
	\]
	
	Hence  $\mathbb{E}[Z_1]=0$.
\end{proof}

\begin{definition}
	Define $M_n:=\max\{Y_i\colon 0\leq i\leq n \}$ and $t_0:=\max\{\log q\colon q\in B \}>0$ and 
	\[
	\tau_r:=\inf\{ i\geq0\colon Y_i\geq rt_0 \}.
	\]
\end{definition}

\begin{lemma}[Ballot theorem]\label{lem: ballot theorem upper bound}
	For $r\geq 1$, 
	\be\label{eq: ballot thm original form}
	\mathbb{P}\big[Y_j>0,j=1,\cdots,n-1,rt_0\leq Y_n <(r+1)t_0\big]\preceq \frac{r}{n^{3/2}}
	\ee
	and 
	\be\label{eq: ballot thm adapted form}
	\mathbb{P}\big[Y_j<t_0,j=1,\cdots,n-1,-(r+1)t_0<Y_n\leq -rt_0\big]\preceq \frac{r}{n^{3/2}}
	\ee
\end{lemma}

For Lemma \ref{lem: ballot theorem upper bound} we need Theorem 8 and 9 from \cite{Addario-Berry2008ballot_thm_published_version}. As in \cite{Addario-Berry2008ballot_thm_published_version} we say a random variable $U$ is non-lattice if there is no real number $\lambda>0$ such that $\lambda U$ is an integer-valued random variable. 
\begin{theorem}[Theorem 8 in \cite{Addario-Berry2008ballot_thm_published_version}]
	\label{thm: ballot thm in the published version}
	Suppose $U$ satisfies $\mathbb{E}[U]=0$, $\textnormal{Var}(U)>0$ and $\mathbb{E}[U^{2+\alpha}]<\infty$ for some $\alpha>0$, and $U$ is non-lattice. Then for any fixed $\beta>0$, given i.i.d.\ random variables $U_1,U_2,\ldots$ distributed as $U$ with associated partial sums $W_i=\sum_{j=1}^{i}U_j$, for all $k$ such that $0\leq k=O(\sqrt{n})$, 
	\[
	\mathbf{P}\big\{ k\leq W_n\leq k+\beta,W_i>0\,\,\forall\,0<i<n  \big\}=\Theta\bigg(\frac{k+1}{n^{3/2}}\bigg).
	\]
\end{theorem}
Theorem 9 from \cite{Addario-Berry2008ballot_thm_published_version} states a corresponding result for the case of $U$ being lattice.

\begin{proof}[Proof of Lemma \ref{lem: ballot theorem upper bound}]
	The inequality \eqref{eq: ballot thm original form} comes directly from Theorem 8 and 9 \cite{Addario-Berry2008ballot_thm_published_version}.  Actually for the upper bound one can drop the assumption  $k=O(\sqrt{n})$ (for instance see Theorem 1 in the arxiv version \cite{addario2008ballot}. The $n^{1/2}$ there was a typo, it should be $n^{3/2}$.) Similarly by applying Theorem 8 and 9 \cite{Addario-Berry2008ballot_thm_published_version} to the partial sums of $-Z_i$'s, one also has that 
	\be\label{eq: minus sum ballot}
		\mathbb{P}\big[Y_j<0,j=1,\cdots,n-1,-(r+1)t_0<Y_n\leq -rt_0\big]\preceq \frac{r}{n^{3/2}}.
	\ee
	
		By Lemma \ref{lem: symmetric random walk on levels induced by the p_h walk}, the vector $\big(Y_1,\ldots, Y_n\big)$ has the same distribution as $\big(Y_2-Y_1,\ldots,Y_n-Y_1,Y_{n+1}-Y_1\big)$ conditioned on $Y_1$. Hence
	\begin{eqnarray*}
	&&\mathbb{P}\big[Y_j<t_0,j=1,\cdots,n-1,-(r+1)t_0<Y_n\leq -rt_0\big]\\
	&=&\mathbb{P}\big[Y_{j+1}-Y_1<t_0,j=1,\cdots,n-1,-(r+1)t_0<Y_{n+1}-Y_1\leq -rt_0\mid Y_1=-t_0\big]\nonumber\\
		&=&\frac{\mathbb{P}\big[Y_1=-t_0,Y_{j+1}<0,j=1,\cdots,n-1,-(r+2)t_0<Y_{n+1}\leq -(r+1)t_0 \big]}{\mathbb{P}[Y_1=-t_0]} \nonumber\\
		&\preceq&\mathbb{P}\big[Y_j<0,j=1,\cdots,n,-(r+2)t_0<Y_{n+1}\leq -(r+1)t_0\big] 
		\stackrel{\eqref{eq: minus sum ballot}}{\preceq}\frac{(r+1)}{n^{3/2}}\preceq \frac{r}{n^{3/2}},
	\end{eqnarray*}
where in the last step we use $r\geq 1$.	
\end{proof}

\begin{lemma}\label{lem: decay of hitting probability of level r}
	For the first hitting times $\tau_r:=\inf\{ i\geq0\colon Y_i\geq rt_0 \}$ one has the following estimate: for all $r\geq1$,
	\be\label{eq: an upper bound on tau_r=k}
	\mathbb{P}[\tau_r=k]\preceq \frac{r}{k^{3/2}}.
	\ee
\end{lemma}
\begin{proof}
	Since the increments $(Z_i)_{i\geq0}$ are a sequence of i.i.d.\ random variables, the vector $(Z_1,\cdots, Z_n)$ has the same distribution as $(Z_n,Z_{n-1},\cdots,Z_1)$. Thus $(Y_1,Y_2,\cdots,Y_n)$ as the partial sum of $(Z_1,\cdots, Z_n)$ has the same distribution as $(Z_n,Z_n+Z_{n-1},\cdots,Z_n+\cdots+Z_1)=(Y_n-Y_{n-1},Y_n-Y_{n-2},\cdots,Y_n-Y_0)$, written as
	\be\label{eq: time reversal for i.i.d. walks}
	(Y_1,Y_2,\cdots,Y_n)\stackrel{\mathscr{D}}{=}(Y_n-Y_{n-1},Y_n-Y_{n-2},\cdots,Y_n-Y_0).
	\ee
	
	Therefore 
	\begin{eqnarray*}
		\mathbb{P}[\tau_r=k]&=&\mathbb{P}[Y_k\geq rt_0,Y_j<rt_0,j=0,1,\cdots,k-1]\nonumber\\
		&\leq &\mathbb{P}\big[Y_k-Y_j>0,j=0,\cdots,k-1,Y_k-Y_0\in[rt_0,(r+1)t_0)\big]\nonumber\\
		&\stackrel{\eqref{eq: time reversal for i.i.d. walks}}{=}&\mathbb{P}\big[Y_j>0,j=1,\cdots,k-1,Y_k\in[rt_0,(r+1)t_0)\big]\nonumber\\
		&\stackrel{\eqref{eq: ballot thm original form}}{\preceq}&\frac{r}{k^{3/2}}.\qedhere
	\end{eqnarray*}
\end{proof}

\begin{lemma}\label{lem: probability of Y_n=0 and M_n=r}
	One has that 
	\be\label{eq: return to level 0}
	\mathbb{P}\big[M_n\in[rt_0,(r+1)t_0), Y_n=0\big]\preceq \frac{(r\vee1)^{3/2}}{n^{3/2}}, 0\leq r\leq \frac{n}{2}.
	\ee
\end{lemma}
\begin{proof}
	We will prove the conclusion for $1\leq r\leq n/2$, the case of $r=0$ being similar and omitted.

	Note that on the event  $\{ M_n\in[rt_0,(r+1)t_0),Y_n=0 \}$, $r\leq \tau_r\leq n-r$ for $1\leq r\leq n/2$.
	Using the strong Markov property of $(Y_n)_{n\geq0}$, by  conditioning on $\tau_r,Y_{\tau_r}$ one has that 
	\be\label{eq: 1.13}
	\mathbb{P}\big[M_n\in[rt_0,(r+1)t_0),Y_n=0\big]
	\leq \sum_{k=r}^{n-r}\mathbb{P}[\tau_r=k]\mathbb{P}\big[ Y_j<t_0,j=1,\cdots,n-k, -(r+1)t_0< Y_{n-k}\leq -rt_0\big].
	\ee

	Therefore 	
	\begin{eqnarray}\label{eq: 2.14}
		\mathbb{P}\big[M_n\in[rt_0,(r+1)t_0),Y_n=0\big]&\stackrel{\eqref{eq: 1.13},\eqref{eq: ballot thm adapted form}}{\leq}&
		\sum_{k=r}^{n-r}\mathbb{P}[\tau_r=k]c\frac{r}{(n-k)^{\frac{3}{2}}}\nonumber\\
		&\stackrel{\eqref{eq: an upper bound on tau_r=k}}{\leq }& 	c_1 \frac{r}{n^{3/2}}\sum_{k=r}^{n/2}\mathbb{P}[\tau_r=k]+\sum_{k=n/2}^{n-r}c_2\frac{r}{k^{3/2}}c_3\frac{r}{(n-k)^{\frac{3}{2}}}\nonumber\\
		&\leq& c_1 \frac{r}{n^{3/2}}+ c_4\frac{r^2}{n^{3/2}}\sum_{k=n/2}^{n-r}\frac{1}{(n-k)^{\frac{3}{2}}}\nonumber\\
		&\leq & c_5\frac{r^{3/2}}{n^{3/2}}.\qedhere
	\end{eqnarray}
\end{proof}

\begin{proof}[Proof of Theorem \ref{thm: main theorem on a_n for amenable, nonuni}]
	
	Write $L_r(o)=\{ v\in V \colon \log \Delta(o,v)\in[rt_0,(r+1)t_0] \}$. 
	Let $x\in V$ be the first vertex in  $L_r(o)$ visited by the $p_h$-walk $(S_i)_{0\leq i\leq n}$. 
		Consider the set $\Gamma_x o$. By Lemma \ref{lem: uni of Haar msr via stabilizer},
	\[
	\frac{|\Gamma_x o|}{|\Gamma_ox|}=\frac{m(x)}{m(o)}=\Delta(o,x)\geq e^{rt_0}.
	\]
	Hence $|\Gamma_xo|\geq|\Gamma_ox|\cdot e^{rt_0}\geq e^{rt_0}$. On the event  $\{M_n\in[rt_0,(r+1)t_0),Y_n=0\}$, by the $\Gamma$-invariance of the $p_h$-walk, the vertices in the set $\Gamma_x o$ are equally likely to be the endpoints of the $p_h$-walk. Hence 
	\[
	\mathbb{P}\big[S_n=o\mid M_n\in[rt_0,(r+1)t_0),Y_n=0 \big]\leq \frac{1}{|\Gamma_xo|}\leq e^{-rt_0}.
	\]
	Therefore  for the $p_h$-walk $(S_n)_{n\geq0}$ starting from $o$,  by Lemma \ref{lem: probability of Y_n=0 and M_n=r} one has that 
	\begin{eqnarray*}
		\mathbb{P}[S_n=o]&\leq&\mathbb{P}\big[M_n\in[0,t_0),Y_n=0\big]\nonumber\\
		&&+\sum_{r=1}^{n/2}\mathbb{P}\big[S_n=o\mid M_n\in[rt_0,(r+1)t_0),Y_n=0 \big]\times \mathbb{P}\big[M_n\in[rt_0,(r+1)t_0),Y_n=0\big]\nonumber\\
		&\leq& c\frac{1}{n^{3/2}}+\sum_{r=1}^{n/2} c_5\frac{r^{3/2}}{n^{3/2}}e^{-rt_0}\preceq \frac{1}{n^{3/2}}.
	\end{eqnarray*}
This establishes $p_h^{(n)}(o,o)\preceq n^{-\frac{3}{2}}$ and then by \eqref{eq: n-step transition prob of p_h is just a_n} we are done.	
\end{proof}

\section{Proof of Theorem \ref{thm: f_n over u_n is at least polynomial decay}}\label{sec: 3}

We begin with some setup and  notation for this section.
\begin{enumerate}
	\item Suppose $G=(V,E)$ is a transitive and $\Gamma$ is a closed, transitive, nonunimodular subgroup of automorphisms. Let $d$ be the degree of $G$. For $x,y\in V(G)$, let $\mathrm{dist}(x,y)$ denote the graph distance between $x$ and $y$ in $G$.
	\item For $o\in V$,  let $\mathscr{L}_n(o)=\{ (v_0,v_1,\cdots,v_n)\colon v_0=v_n=o,v_i\sim v_{i+1} \textnormal{ for }i=0,\ldots,n-1 \}$ be the set of cycles \notion{rooted at $o$} with length $n$ in $G$.
	
	
	\item Let $(X_n)_{n\geq0}$ be a simple random walk on $G$. Denote by $\mathbb{P}_o$ the law of $(X_n)_{n\geq0}$ when the walk starts from $X_0=o$. For $w\in\mathscr{L}_n(o)$,  $\mathbb{P}_o[(X_0,\cdots,X_n)=w]=\frac{1}{d^n}$ is the probability of traveling along the particular path  $w$ by a simple random walk for the first $n$ steps.

\end{enumerate}

\begin{definition}
	From item 3 in the above, conditioned on $X_n=X_0=o$, the trajectory $(X_0,\ldots,X_n)$ can be sampled from  $\mathscr{L}_n(o)$ uniformly at random. We  denote the law of the conditional trajectory by $\mathbf{P}_{n,o}$. Let $\mathbf{E}_{n,o}$  denote the corresponding expectation. 
\end{definition}

\begin{lemma}\label{lem: 3.3}
For a transitive graph $G=(V,E)$ and $o,x\in V$,  one has that
	\[
	\mathbf{P}_{n,o}[x\in w]=\mathbf{P}_{n,x}[o\in w].
    \]
\end{lemma}
The proof of Lemma \ref{lem: 3.3} is a routine application of  the reversibility and symmetry of the random walk together with the transitivity of the graph. Hence the proof is omitted.

	For $k\in\mathbb{Z}$, define that $L_k(x):=\big\{ y\in V\colon  \log \Delta(x,y)\in[kt_0,(k+1)t_0] \big \}$. 

\begin{proposition}\label{prop: TMTP for the first moment intersect with a level}
	For $0\leq k\leq n$ one has that
	\be\label{eq: 3.9}
	\mathbf{E}_{n,x}\big[ |w\cap L_k(x)| \big]\leq ne^{-t_0k},
	\ee
	where $|w\cap L_k(x)|$ is the number of vertices in the intersection of $w$ with $L_k(x)$.
\end{proposition}
\begin{proof}

	Define a function $f:V\times V\to[0,\infty)$ by
	\[
	f(x,y)
	=\mathbf{1}_{y\in L_k(x)}\cdot \mathbf{P}_{n,x}[y\in w]=\mathbf{1}_{y\in L_k(x)}\cdot \mathbf{E}_{n,x}\big[\mathbf{1}_{y\in w}  \big].
	\]
	The function  $f$ is $\Gamma$-diagonally invariant by the $\Gamma$-diagonal invariance of  the modular function $\Delta$ (Lemma \ref{lem: properties of modular func}) and transitivity of the graph. By the mass-transport principle, we have
	\begin{eqnarray}
		\mathbf{E}_{n,x}\big[ |w\cap L_k(x)| \big]
		&=& \sum_{y\in V} \mathbf{1}_{y\in L_k(x)}\mathbf{E}_{n,x}\big[   \mathbf{1}_{y\in w}\big]=\sum_{y\in V}f(x,y)\nonumber\\
		&=& \sum_{y\in V}f(y,x)\Delta(x,y)= \sum_{y\in V} \mathbf{1}_{x\in L_k(y)}
		\mathbf{P}_{n,y}[x\in w] \cdot \Delta(x,y)
	\end{eqnarray}	
If  $\mathbf{1}_{x\in L_k(y)}=1$, then  $\log\Delta(y,x)\in[kt_0,(k+1)t_0]$ and $\log\Delta(x,y)=-\log\Delta(y,x)\in \big[-(k+1)t_0,-kt_0\big]$.  This implies that if $\mathbf{1}_{x\in L_k(y)}=1$, then $y\in L_{-k-1}(x)$ and $\Delta(x,y)\leq e^{-kt_0}$. Therefore 
\be
\mathbf{E}_{n,x}\big[ |w\cap L_k(x)| \big]\leq e^{-kt_0}\sum_{y\in V} \mathbf{1}_{y\in  L_{-k-1}(x)} \cdot \mathbf{P}_{n,y}[  x\in w]
\ee
By Lemma \ref{lem: 3.3} one has that
\be\label{eq: 3.5}
	\mathbf{E}_{n,x}\big[ |w\cap L_k(x)| \big]
	\leq e^{-kt_0}\sum_{y\in V} \mathbf{1}_{y\in L_{-k-1}(x)} \cdot \mathbf{P}_{n,x}[  y\in w].
\ee
	
	Since $w_0=w_n=x$ for all $w\in \mathscr{L}_{n}(x)$, one has $|w|\leq n$ and thus \eqref{eq: 3.9}:
	\[
	\mathbf{E}_{n,x}\big[ |w\cap L_k(x)| \big]\stackrel{\eqref{eq: 3.5}}{\leq} e^{-kt_0}\mathbf{E}_{n,x}\big[ \bigl|w\cap \big(L_{-k-1}(x)\big)\bigr|  \big]
	\leq e^{-kt_0}\mathbf{E}_{n,x}\big[ |w| \big]\leq  ne^{-t_0k}. \qedhere
	\]
\end{proof}

\begin{lemma}\label{lem: comparison of u_2s+t and u_2s}
	Suppose $G=(V,E)$ is a transitive graph with spectral radius $\rho$.
	There is a constant $c_1=c_1(G)>0$ and $n_1=n_1(G)\geq 0$ such that for all $k\geq 1$ and $n\geq2k+n_1$,
	\be\label{eq: 3.13}
	u_{n-2k}\geq c_1u_{n}\rho^{-2k}
	\ee
	
\end{lemma}
\begin{proof}
	We first review a classical application of Cauchy--Schwarz inequality from the proof Lemma 10.1 in \cite{woess2000random}. Let $(\cdot,\cdot)$ denote the standard inner product on $l^2(V)$. Let $f:V\to\mathbb{R}$ be a non-negative function with finite support. Let $P$ be the transition operator associated with simple random walk on $G$. Then $P$ is a self-adjoint operator on $l^2(V)$ and  $(P^nf,P^nf)=(f,P^{2n}f)$ is finite for each $n$. Using Cauchy--Schwarz inequality one has that 
	\[
	(P^{n+1}f,P^{n+1}f)^2=(P^nf,P^{n+1}f)^2\leq (P^nf,P^nf)(P^{n+2}f,P^{n+2}f).
	\]
	Hence the sequence $\frac{(P^{n+1}f,P^{n+1}f)}{(P^nf,P^nf)}$ is increasing. The limit is then equal to $(P^nf,P^nf)^{1/n}$. Hence by taking $f=\mathbf{1}_{x}$ one has that
	\be\label{eq: Cauchy Schwarz inequality application}
	u_2\leq \frac{u_{2k+2}}{u_{2k}}\leq \lim_{n\to \infty}u_{2n}^{1/n}=\rho^2. 
	\ee

	Now if $n$ is even, then \eqref{eq: 3.13} actually holds for $c_1=1$: for $n=2m$, by \eqref{eq: Cauchy Schwarz inequality application} one has that
	\[
	\frac{u_{n}}{u_{n-2k}}=\frac{u_{2m}}{u_{2m-2k}}=\prod_{j=1}^{k}\frac{u_{2m-2k+2j}}{u_{2m-2k+2(j-1)}}\leq \rho^{2k}.
	\]

	Second, if $n$ is odd, say $n=2m+1$,  then we can assume that there exists a smallest odd number $2l+1>0$ such that $u_{2l+1}>0$; otherwise \eqref{eq: 3.13} is trivial because both sides are zero.  In particular, $u_{2l+1}=\sum_{j=1}^{2l+1}f_ju_{2l+1-j}=f_{2l+1}$ (when $j<2l+1$, if $j$ is odd, then $f_j\leq u_j=0$; if $j$ is even, then  $u_{2l+1-j}=0$). 
	
	The inequality  (2.10) in Lemma 1 of \cite{BIK2007} says that 
	\be\label{eq: 3.14}
	u_n=u_{2m+1}\leq u_{2m}. 
	\ee

	Take $n_1=2l+1$. For $k\geq 1$ and $n-2k=2m+1-2k\geq n_1=2l+1$, one has that 
	\begin{eqnarray}
		u_{n-2k}&=&u_{2m+1-2k}\geq  f_{2l+1}u_{2m-2k-2l}=u_{2l+1}u_{2m-2k-2l}\nonumber\\
		&\geq& u_{2l+1} u_{2m}\rho^{-(2k+2l)}\nonumber\\
		&\stackrel{\eqref{eq: 3.14}}{\geq}&  \frac{u_{2l+1}}{\rho^{2l}}u_{2m+1}\rho^{-2k}=\frac{u_{2l+1}}{\rho^{2l}}u_{n}\rho^{-2k},
	\end{eqnarray}
	where in the second step we use \eqref{eq: 3.13} with $c_1=1$ for the even case that we already proved.
	
	Since  $\frac{u_{2l+1}}{\rho^{2l}}\stackrel{\eqref{eq: 3.14}}{\leq}\frac{u_{2l}}{\rho^{2l}}\leq 1$, taking $c_1=\frac{u_{2l+1}}{\rho^{2l}}$ and $n_1=2l+1$ we have the desired conclusion for odd $n$. 
\end{proof}

\begin{proof}[Proof of Theorem \ref{thm: f_n over u_n is at least polynomial decay}]
	Write $H_k^+(x)=\bigcup_{n\geq k}L_k(x)$. By the definition of $t_0$, a path $w\in \mathscr{L}_n(x)$ intersects with $H_k^+(x)$ if and only if $|w\cap L_{k}(x)|\geq 1$. 
	By Proposition \ref{prop: TMTP for the first moment intersect with a level}, for $k\geq0$,
	\be\label{eq: exponential decay of prob intersecting with a high level}
	\mathbf{P}_{n,x}[w\cap H^+_{k}(x)\neq \emptyset]=\mathbf{P}_{n,x}[|w\cap L_{k}(x)|\geq 1]\leq \mathbf{E}_{n,x}[|w\cap L_{k}(x)|]\leq ne^{-kt_0}
	\ee
	Take $C>0$ large such that $ne^{-kt_0}\leq \frac{1}{2}$ for $k\geq k(n):=\lfloor C\log n\rfloor$. 
	Hence 
	\be\label{eq: 3.18}
	\mathbf{P}_{n,x}[w\cap H^+_{k}(x)=\emptyset]\geq \frac{1}{2},\,\,\,\forall\,\,k\geq k(n).
	\ee
	
	Note that there is a path $\gamma$ of length $k=k(n)$ from $x$ to some $y$ such that $\mathrm{dist}(y,x)=k+1$ and $\log \Delta(x,y)=-(k+1)t_0$.
	
	Suppose we first travel from $x$ to $y$ along $\gamma$ in the first $k+1$ steps, and then do an excursion from $y$ to $y$ such that the cycle has length $n-2(k+1)$ and doesn't intersect with $L_k(y)$, then travel from $y$ to $x$ along the reversal of $\gamma$. Then  we come back to $x$ for the first time at time $n$. Recall that we denote by $f_n$ the first return probability for simple random walk. 
	Since $\lfloor C\log n\rfloor\geq \lfloor C\log \big(n-\lfloor C\log n\rfloor\big)\rfloor$,  for all $n$ sufficiently large the first return probability $f_n$ satisfies
	\begin{eqnarray*}
		f_n&\geq& \frac{1}{d^{k(n)+1}}\times u_{n-2k(n)-2}\times \mathbf{P}_{n-2k(n)-2,y}[w\cap L_{k(n)}(y)=\emptyset]\times \frac{1}{d^{k(n)+1}}\nonumber\\
		&\stackrel{\eqref{eq: 3.18}}{\geq } & \frac{1}{d^{2k(n)+2}}\times u_{n-2k(n)-2}\times \frac{1}{2}\nonumber\\
		&\stackrel{\mathrm{Lem.}\ref{lem: comparison of u_2s+t and u_2s}}{\geq} & \frac{c_1}{2(d\rho)^{2k(n)+2}}u_n\succeq \frac{1}{n^c}u_n.
	\end{eqnarray*}
It is easy to see that $\rho\geq \frac{1}{d}$; but the equality can't hold in our case. In fact Theorem 6.10 in \cite{LP2016} says that for a connected, regular, infinite graph with spectral radius $\rho$ and degree $d$, one always has that $\rho\cdot d \geq 2\sqrt{d-1}>1$.
\end{proof}

\section{Extensions to quasi-transitive graphs}\label{sec: quasi-transitive case}

Suppose $\Gamma\subset \mathrm{Aut}(G)$ is a closed subgroup of automorphisms of a locally finite, connected graph $G=(V,E)$. For $v\in V$, let $\Gamma v=\{\gamma v\colon \gamma\in\Gamma\}$ denote the orbit of $v$ under $\Gamma$. Let $G/\Gamma=\{\Gamma v\colon v\in V\}$ be the set of orbits for the action of $\Gamma$ on $G$. We say $\Gamma$ is \textbf{quasi-transitive} if $G/\Gamma$ is a finite set. 
In this section we extend Theorem \ref{thm: main theorem on a_n for amenable, nonuni} and \ref{thm: f_n over u_n is at least polynomial decay} to the quasi-transitive case.

\subsection{Extension of Theorem \ref{thm: main theorem on a_n for amenable, nonuni}}
Recall that for a graph $G=(V,E)$ with spectral radius $\rho$ and a vertex $x\in V$, we denote by $u_n(x)=\mathbb{P}[X_n=x\mid X_0=x]$ the $n$-step return probability for simple random walk $(X_i)_{i\geq0}$ on $G$ and $a_n(x)=\frac{u_n(x)}{\rho^n}$. 
\begin{theorem}\label{thm: main theorem on a_n for amenable, nonuni, quasi-transitive}
	Suppose $G=(V,E)$ is a locally finite, connected graph with  a closed, quasi-transitive, amenable and nonunimodular subgroup of automorphisms. Then $a_n(x)\preceq n^{-\frac{3}{2}}$, $\forall \,x\in V$.
\end{theorem}
The idea for the quasi-transitive case is the same as the transitive case: find a $\rho$-harmonic function $h$ and then consider the associated $p_h$-walk.  

\subsubsection{The $\rho$-harmonic function $h$ in the quasi-transitive case}
We first set up some notation. Throughout this subsection we assume $G=(V,E)$ is a connected, infinite graph with
$\Gamma$ being a closed, amenable, quasi-transitive subgroup of $\mathrm{Aut}(G)$. Let $\mathcal{O}=\{o_1,\ldots,o_L\}$ be a complete set of representatives in $V$ for the orbits of $\Gamma$. Let $I=\{1,\ldots, L\}$ be the index set. For $x\in V$, let $d_x$ be the degree of $x$. We also write $d_x=d_i$ when $x\in \Gamma o_i$ since the degrees of the vertices in the same orbit are the same.  Recall that $m(x)=|\Gamma_x|$ is the left-Haar measure of the stabilizer $\Gamma_x$.
Recall that in the case $\Gamma$ acts transitively on $G$ we use Theorem 1(b) from \cite{soardi1990amenability} to establish the $\rho$-harmonicity of the associated function $h$. Here we need a natural extension of Theorem 1(b) from \cite{soardi1990amenability}, namely Theorem 1(b) from \cite{Saloff-Coste_and_Woess1996norm_transition}. 

Let $A=(a(i,j))_{i,j\in I}$ be the matrix as defined in \cite{Saloff-Coste_and_Woess1996norm_transition}, namely, 
\[
a(i,j)=\sum_{y\in\Gamma o_j}\frac{\mathbf{1}_{\{ y\sim x \}}}{d_x}\sqrt{\frac{d_x}{d_y} \frac{m(y)}{m(x)} }
=\sum_{y\in\Gamma o_j,y\sim x}\frac{1}{\sqrt{d_xd_y}}\sqrt{\frac{m(y)}{m(x)}}, x\in\Gamma o_i,\,i,j\in I.
\]
Note that the $\Gamma$-invariance of $\Delta(x,y)=\frac{m(y)}{m(x)}$ (Lemma \ref{lem: properties of modular func}) yields that $a(i,j)$ does not depend on the choice of $x\in \Gamma o_i$. 
Obviously the matrix $A$ is  irreducible and nonnegative. Hence by Perron-Frobenius theorem there is a positive  vector $\vec{v}=(v_1,\ldots,v_L)^T$ associated to the largest eigenvalue $\rho(A)$ (we normalize $\vec{v}$ to have $l_2$-norm $1$). Theorem 1(b) of \cite{Saloff-Coste_and_Woess1996norm_transition} says that for amenable $\Gamma$ one has that $\rho=\lambda(A)$, where $\rho$ is the spectral radius of $G$ and $\lambda(A)$ is the largest eigenvalue of the finite matrix $A$. Hence for all $i\in I$ the eigenvalue equation becomes
\be\label{eq: eigenvector equation}
\rho v_i=\sum_{j=1}^{L}a(i,j)v_j.
\ee
Also Lemma 3(1) of \cite{Saloff-Coste_and_Woess1996norm_transition} says that the matrix $A$ is symmetric, i.e.,
\be\label{eq: symmetry of the matrix A}
a(i,j)=a(j,i).
\ee
\begin{definition}\label{def: quasi-transitive rho-harmonic function h}
	Define $v:V\to (0,\infty)$ by setting  $v(x)=v_i$ for $x\in\Gamma o_i$. Let $h:V\to(0,\infty)$ be given by $h(x)=v(x)\sqrt{\frac{m(x)}{d_x}}$.
\end{definition}
We will see that this function $h$ is $\rho$-harmonic and hence as before one can define the associated $p_h$-walk on $G$ via the transition probabilities:
\[
p_h(x,y)=\frac{p(x,y)h(y)}{\rho h(x)}=\frac{\mathbf{1}_{y\sim x}}{d_x}\cdot \frac{h(y)}{\rho h(x)}.
\]

\begin{proposition}\label{prop: rho harmonicity in the quasitransitive case}
	The function $h$ defined in Definition \ref{def: quasi-transitive rho-harmonic function h} is $\rho$-harmonic. The associated $p_h$-walk on $G$ is reversible with respect to $v^2m$.  
\end{proposition}
\begin{proof}
We first verify that $h$ is $\rho$-harmonic. For an arbitrary vertex $x\in V$, say $x\in \Gamma o_i$, we have that
\begin{eqnarray*}
	\frac{1}{d_x}\sum_{ y\sim x}h(y)&=&\sum_{j=1}^{L}\sum_{ y\in \Gamma o_j,y\sim x}\frac{1}{d_x}h(y)=\sum_{j=1}^{L}\sum_{ y\in \Gamma o_j,y\sim x}\frac{1}{d_x}v_j\sqrt{\frac{m(y)}{d_y}}\\
	&=&\sum_{j=1}^{L}\frac{\sqrt{m(x)}}{\sqrt{d_x}}v_j\sum_{y\in\Gamma o_j,y\sim x}\frac{1}{\sqrt{d_xd_y}}\sqrt{\frac{m(y)}{m(x)}}\\
	&=&\frac{\sqrt{m(x)}}{\sqrt{d_x}}\sum_{j=1}^{L}v_ja(i,j)\stackrel{\eqref{eq: eigenvector equation}}{=}\frac{\sqrt{m(x)}}{\sqrt{d_x}}\cdot \rho v_i=\rho\cdot h(x).
\end{eqnarray*}

For reversibility of the $p_h$-walk, it is also straightforward to verify that $v(x)^2m(x)p_h(x,y)=v(y)^2m(y) p_h(y,x)$ and  details are skipped.
\end{proof}

\subsubsection{The increments have mean zero when starting from the stationary distribution.}

Let $(S_n)_{n\geq0}$ be a $p_h$-walk on $G$ and let $Y_n=\log\Delta(S_0,S_n)$ be the associated process on $\mathbb{R}$. Here we recall  the modular function $\Delta(x,y)=\frac{|\Gamma_yx|}{|\Gamma_xy|}$ given in Definition \ref{def: modular function}.  Using Lemma \ref{lem: uni of Haar msr via stabilizer} one has that $\Delta(x,y)=\frac{m(y)}{m(x)}$. Hence the cocycle identity $\Delta(x,y)\Delta(y,z)=\Delta(x,z)$ still holds in the quasi-transitive case. By the cocycle identity the increment at time $i$ of the process $(Y_n)_{n\geq0}$ is   $\log \Delta(S_i,S_{i+1})$. Note that the distribution of this increment $\log \Delta(S_i,S_{i+1})$ depends (and only depends) on the orbit of $S_i$. So in order to have mean-zero increments in the long run one must have mean-zero increments when starting from the stationary distribution and this is indeed the case (Prop.\,\ref{prop: zero mean increment for stationary distributed initial starting point}). 

\begin{definition}\label{def: stationary distribution for the quasitransitive p_h walk}
	Recall that the vector $\vec{v}=(v_1,\ldots,v_L)^T$ with $l_2$-norm $1$ is the unique eigenvector associated with the largest eigenvalue $\rho$ of the matrix $A$.
	Define $\pi=(\pi_i)_{i\in I}$ by  $\pi_i=v_i^2,i\in I$.  
\end{definition}

\begin{proposition}\label{prop: zero mean increment for stationary distributed initial starting point}
	 The measure $\pi=(\pi_i)_{i\in I}$ is the stationary probability measure for the induced chain on $I$ of the $p_h$-walk. Let $(S_n)_{n\geq0}$ be a $p_h$-walk on $G$ with starting point $S_0$ sampled from the measure $\pi$. Then
	\be\label{eq: mean zero when start with a stationary vertex}
	\mathbb{E}\big[ \log 	\Delta(S_0,S_1) \big]=0
	\ee
	

\end{proposition}

The following lemma is an analogue of Lemma \ref{lem: ratio of m(x) for neighboring vertices}.
\begin{lemma}
	Write $B_{i,j}=\big\{ \frac{m(y)}{m(x)} \colon x\in\Gamma o_i, y\in \Gamma o_j,x\sim y \big\}$. For $q\in B_{i,j}$, let $N_{i,j,q}=\{ y\colon y\in \Gamma o_j, y\sim o_i, \frac{m(y)}{m(o_i)}=q \}$. Then
	 $q\in B_{i,j}$ if and only if $q^{-1}\in B_{j,i}$ and  
	\be\label{eq: relation on the number of neighboring points with a fixed ratio}
	\# N_{i,j,q}=\frac{1}{q} \# N_{j,i,q^{-1}}.
	\ee
\end{lemma}
\begin{proof}
	Let  $f: V\times V\to[0,\infty)$  be the indicator function given by 
	\[
	f(x,y)=\mathbf{1}_{\{ x\in\Gamma o_i,\, y\in \Gamma o_j,\,x\sim y,\, \frac{m(y)}{m(x)}=q  \}}.
	\]
	Obviously $f$ is $\Gamma$-diagonally invariant. Hence by the mass-transport principle (Prop.\,\ref{prop: TMTP}) one has \eqref{eq: relation on the number of neighboring points with a fixed ratio}:
	\[
	\# N_{i,j,q}=\sum_{z\in \Gamma o_j}f(o_i,z)=\sum_{y\in \Gamma o_i}f(y,o_j)\frac{m(y)}{m(o_j)}=\frac{1}{q}\# N_{j,i,q^{-1}}.\qedhere
	\]
\end{proof}

\begin{proof}[Proof of Proposition \ref{prop: zero mean increment for stationary distributed initial starting point}]
	The $p_h$-walk $(S_n)_{n\geq0}$ induces a Markov chain on the index set $I$ with transition probability $\widetilde{P}_h$ given by: $\forall\,i,j\in I,$
	\begin{eqnarray*}
	\widetilde{P}_h(i,j)&=&\mathbb{P}[S_1\in\Gamma o_j\mid S_0\in \Gamma o_i]=\sum_{y\in\Gamma o_j,y\sim o_i} p_h(o_i,y)\nonumber\\
	&=&\sum_{y\in\Gamma o_j,y\sim o_i}\frac{h(y)}{d_i\cdot \rho\cdot h(o_i)}=\sum_{y\in\Gamma o_j,y\sim o_i} \frac{1}{\rho}\cdot \frac{v_j}{v_i}\cdot \frac{1}{\sqrt{d_id_j}}\cdot \sqrt{\frac{m(y)}{m(o_i)}}\nonumber\\
	&=&\sum_{q\in B_{i,j}} \frac{1}{\rho}\cdot \frac{v_j}{v_i}\cdot \frac{1}{\sqrt{d_id_j}}\cdot \sqrt{q} \cdot \# N_{i,j,q}. 
	\end{eqnarray*}
	Now we verify the stationarity of $\pi$ for $\widetilde{P}_h$:
	\begin{eqnarray*}
		\sum_{i=1}^{L}\pi_i\widetilde{P}_h(i,j)&=&\sum_{i=1}^{L}v_i^2\sum_{q\in B_{i,j}}\frac{1}{\rho}\cdot \frac{v_j}{v_i}\cdot \frac{1}{\sqrt{d_id_j}}\cdot \sqrt{q} \cdot \# N_{i,j,q}\nonumber\\
		&=& v_j\sum_{i=1}^{L}v_i\sum_{q^{-1}\in B_{j,i}} \frac{1}{\rho}\cdot \frac{1}{\sqrt{d_id_j}}\cdot \sqrt{q} \cdot \# N_{i,j,q}\nonumber\\
		&\stackrel{\eqref{eq: relation on the number of neighboring points with a fixed ratio}}{=}&v_j \sum_{i=1}^{L}v_i\sum_{q^{-1}\in B_{j,i}} \frac{1}{\rho}\cdot \frac{1}{\sqrt{d_id_j}}\cdot \sqrt{q^{-1}}\cdot \# N_{j,i,q^{-1}}\nonumber\\
		&=&v_j\sum_{i=1}^{L}v_ia(j,i)\cdot \frac{1}{\rho} \stackrel{\eqref{eq: symmetry of the matrix A}}{=} v_j\sum_{i=1}^{L}v_ia(i,j)\cdot \frac{1}{\rho}\stackrel{\eqref{eq: eigenvector equation}}{=}v_j^2=\pi_j.
	\end{eqnarray*}

Finally we verify \eqref{eq: mean zero when start with a stationary vertex}:
\begin{eqnarray*}
		\mathbb{E}\big[ \log 	\Delta(S_0,S_1) \big]
		&=&\sum_{i=1}^{L}\pi_i \sum_{j=1}^{L}\sum_{q\in B_{i,j}} \sum_{y\in N_{i,j,q}}p_h(o_i,y) \log q\nonumber\\
		&=&\sum_{i=1}^{L}v_i^2 \sum_{j=1}^{L}\sum_{q\in B_{i,j}} \sum_{y\in N_{i,j,q}} \Big(\frac{1}{\rho}\cdot \frac{v_j}{v_i}\cdot \frac{1}{\sqrt{d_id_j}}\cdot \sqrt{q}\Big) \cdot \log q\nonumber\\
		&=&\sum_{i=1}^{L}v_i^2 \sum_{j=1}^{L}\sum_{q\in B_{i,j}} \frac{\#N_{i,j,q}}{\rho}\cdot \frac{v_j}{v_i}\cdot \frac{1}{\sqrt{d_id_j}}\cdot \sqrt{q} \log q\nonumber\\
		&=&\frac{1}{\rho}\sum_{i=1}^{L} \sum_{j=1}^{L}\sum_{q\in B_{i,j}} \big(\sqrt{q}\#N_{i,j,q}\big)\cdot (v_iv_j)\cdot \frac{1}{\sqrt{d_id_j}}\cdot \log q.
\end{eqnarray*}
By \eqref{eq: relation on the number of neighboring points with a fixed ratio}, the term $t_{i,j,q}=\big(\sqrt{q}\#N_{i,j,q}\big)\cdot (v_iv_j)\cdot \frac{1}{\sqrt{d_id_j}}\cdot \log q$ exactly cancels the term $t_{j,i,q^{-1}}$ and one obtains \eqref{eq: mean zero when start with a stationary vertex}.
\end{proof}

\subsubsection{Proof of Theorem \ref{thm: main theorem on a_n for amenable, nonuni, quasi-transitive}}
Similar to the transitive case, to prove Theorem \ref{thm: main theorem on a_n for amenable, nonuni, quasi-transitive} it suffices to show the following analogue of Lemma \ref{lem: probability of Y_n=0 and M_n=r}. 
\begin{lemma}\label{lem: prob of M_n=rt_0, Y_n=0 for quasi-transitive}
	Let $(S_n)_{n\geq0}$ be a $p_h$-walk on $G$ starting from a random point in $\mathcal{O}$ sampled according to the measure $\pi$ from Definition \ref{def: stationary distribution for the quasitransitive p_h walk}. As the transitive case, let $M_n=\max\{Y_i\colon 0\leq i\leq n \}$ and $t_0=\max\{\log \Delta(x,y) \colon x\sim y  \}>0$.  Let $Y_k=\log \Delta(S_0,S_k),k\geq0$.   
	
	Then 
	\be\label{eq: return to level 0 with maximum in level r for quasi-transitive}
	\mathbb{P}\big[M_n\in[rt_0,(r+1)t_0), Y_n=0\big]\preceq \frac{(r\vee1)^{3/2}}{n^{3/2}}, 0\leq r\leq \frac{n}{2}.
	\ee
\end{lemma}
The proof of Lemma \ref{lem: prob of M_n=rt_0, Y_n=0 for quasi-transitive} follows a similar  strategy for Lemma \ref{lem: probability of Y_n=0 and M_n=r}. Hence we put the details of the proof in the appendix for completeness.

\subsection{Extension of Theorem \ref{thm: f_n over u_n is at least polynomial decay}}

\begin{theorem}\label{thm: main theorem for radii of convergence of F(z) for nonuni, quasi-transitive}
	Suppose $G=(V,E)$ is a locally finite, connected graph with  spectral radius $\rho$  and a closed, quasi-transitive and nonunimodular subgroup of automorphisms. Write $f_n(x)=\mathbb{P}_x[X_n=x,X_i\neq x, \forall \, 1\leq i\leq n-1]$ for the first return probability, where $(X_n)_{n\geq0}$ is a simple random walk on $G$ starting from $X_0=x$. Then 
 there exists $c>0$ such that for all $x\in V, n>0$,
	\[
	f_n(x)\succeq \frac{1}{n^c}u_n(x).
	\]
\end{theorem}
\begin{proof}
	For the quasi-transitive case,  it is easy to see that there are constants $C>1$ and $n_0>0$ such that for any $n\geq n_0$,
	\[
	\frac{1}{C}u_n(y)\leq u_n(x)\leq Cu_n(y),\,\,\forall\,x,y\in V. 
	\]
	Lemma \ref{lem: 3.3} obviously holds for $o,x$ in the same orbit and then similar to Proposition \ref{prop: TMTP for the first moment intersect with a level}  one has that 
	\be\label{eq: intersection within same orbit}
		\mathbf{E}_{n,o}\big[ |w\cap L_k(o)\cap \Gamma o| \big]\leq ne^{-t_0k},
	\ee
	
	 By quasi-transitivity and connectedness of $G$, there is a constant $D>0$ such that for any $o,x\in V$, there is a point $x'=x(o)\in\Gamma o$   such that $\mathrm{dist}(x,x')\leq D$. 
Note that
	\[
\mathbb{P}_o[X_n=o,X_i=x\textnormal{ for some }i<n]\preceq \mathbb{P}_o[X_{n+2\mathrm{dist}(x,x')}=o,X_j=x'\textnormal{ for some }j<n+2\mathrm{dist}(x,x')].
	\]
Since  by Lemma \ref{lem: comparison of u_2s+t and u_2s} $u_n(o)\asymp u_{n+2t}(o)$ for $t\leq D$,  one has that
\[
\mathbf{P}_{n,o}[x\in w]\preceq \mathbf{P}_{n+2\mathrm{dist}(x,x'),o}[x'\in w].
\]
Summing this over $x\in L_k(o)$ (the corresponding $x'\in L_{k+t'}(o)$ for some $t'$ satisfies $|t'|\leq D$, and each $x'$ can added up to $|B(o,D)|$ times)
 and using \eqref{eq: intersection within same orbit} one has that 
	\[
		\mathbf{E}_{n,o}\big[ |w\cap L_k(o)| \big]\preceq \big( n+2D\big)e^{-t_0k}\preceq ne^{-t_0k}.
	\]
	The rest is the same as the transitive case. 
\end{proof}

\section{Some nonunimodular examples for Conjecture \ref{conj: a_n is at most n^(-3/2)}}\label{sec: examples}\label{sec: 5 examples for conjecture 1}
We have seen that Conjecture \ref{conj: a_n is at most n^(-3/2)} holds for all transient, transitive,  amenable graphs. In this section we give some nonunimodular examples for which Theorem \ref{thm: main theorem on a_n for amenable, nonuni} applies. Among the following examples, the result $a_n\asymp n^{-3/2}$ for grandparent graphs in Example \ref{example: grandparent graph} and $a_n\preceq n^{-3/2}$ for certain Cartesian products as in Example \ref{example: Cartesian product} seem to be new.

We first recall a simple criterion for the amenability of a subgroup of automorphisms of certain graphs. If $G$ has infinitely many ends, the following  proposition from \cite{soardi1990amenability} gives a way to determine amenability of a closed transitive subgroup of automorphisms. 
\begin{proposition}[Proposition 2 of \cite{soardi1990amenability}]\label{prop: amenability equiv fix end}
	Let $\Gamma$ be a closed, transitive subgroup of $\mathrm{Aut}(G)$ for the graph $G$ and $G$ has infinitely many ends. Then $\Gamma$ is amenable iff it fixes a unique end. 
\end{proposition}

\begin{example}[Toy model]\label{example: toy model}
	Consider a regular tree $\mathbb{T}_{b+1}$ with degree $b+1\geq3$. Let $\xi$ be an end of the tree and $\Gamma_\xi$ be the subgroup of automorphisms that fixes the end $\xi$. Then $\Gamma_\xi$ is  a closed, amenable, nonunimodular, transitive subgroup of $\mathrm{Aut}(\mathbb{T}_{b+1})$. The transitivity can be easily verified. The amenability follows from Proposition \ref{prop: amenability equiv fix end}. The nonunimodularity follows from a simple application of Proposition \ref{prop: unimodularity iff modular func is constant 1}.

\end{example}

Typical examples on nonunimodular transitive graphs are grandparent graphs and Diestel-Leader graphs which we now briefly recall.

\begin{example}[Grandparent graph]\label{example: grandparent graph}
	Let $\xi$ be a fixed end of a regular tree $\mathbb{T}_{b+1}$ ($b\geq2$) as in the toy model. For a vertex $v\in \mathbb{T}_{b+1}$, there is a unique ray $\xi_v:=\{v_0,v_1,v_2,\ldots\}$ representing $\xi$ started at $v_0=v$. Call $v_2$ in the ray $\xi_v$ the ($\xi$-)grandparent of $v$. Add edges between all vertices and their grandparents and the graph $G$ obtained is called a grandparent graph. It is easy to see that $\mathrm{Aut}(G)=\Gamma_\xi$, the subgroup of $\mathrm{Aut}(\mathbb{T}_{b+1})$ that fixes the end $\xi$. Hence Theorem \ref{thm: main theorem on a_n for amenable, nonuni} applies to grandparent graphs. 
	
	In fact from the proof of Theorem \ref{thm: main theorem on a_n for amenable, nonuni} and the underlying tree-like structure of $G$, one has that $a_n\asymp n^{-3/2}$. (The lower bound can be showed by considering the probability that the $p_h$-walk started from $x$ returns to the level $L_0(x)$ at time $n$ without using any vertex of in $L_k(x),k>0$; the tree-like structure then force the returning point at time $n$  in the level $L_0(x)$ must be  $x$ itself.)
\end{example}

\begin{remark}
	For the toy model in Example \ref{example: toy model} or the grandparent graph in Example \ref{example: grandparent graph}, using Lemma \ref{lem: probability of Y_n=0 and M_n=r} and the fact that $a_{2n}\asymp n^{-3/2}$ for these two cases one actually improves \eqref{eq: exponential decay of prob intersecting with a high level} to 
	\[
	\mathbf{P}_{n,x}[w\cap H^+_{k}(x)\neq \emptyset]\preceq (k\vee 1)^{3/2}e^{-kt_0}.
	\]
%
	Hence  there is a large constant $k$ (independent of $n$) such that 
	\be\label{eq: remark for possible refinement}
	\mathbf{P}_{n,x}[w\cap H^+_{k}(x)\neq \emptyset]\leq\frac{1}{2}.
	\ee
	Then one can adapt the proof of Theorem \ref{thm: f_n over u_n is at least polynomial decay} to show that  $f_n\asymp u_n$ (a path with  length of constant order instead of  $\log n$ would suffice). So a natural question is whether such an inequality \eqref{eq: remark for possible refinement} holds for  general nonunimodular, transitive graphs. 
\end{remark}

\begin{example}[Diestel-Leader graph]\label{example: DL graph}
	Woess asked whether there is a vertex-transitive graph that is not roughly isometric to any Cayley graph.
	Diestel and Leader \cite{Diestel_and_Leader2001} construct a family of graphs $DL(q,r)$ and conjectured these graphs are  not roughly isometric to any Cayley graph when $q\neq r$. Later it was proved this is indeed the case \cite{Eskin_Fisher_Whyte2012not_quasi-isometric_Cayley}. Now these graphs are called Diestel-Leader graphs. 
	
	Let $G_1=\mathbb{T}_{q+1},G_2=\mathbb{T}_{r+1}$ be two regular trees with degree $q+1,r+1\geq 3$ respectively. Let $\xi_i$ be an end of $G_i$, $i=1,2$. Let $\Gamma_i$ be the subgroup of $\mathrm{Aut}(G_i)$ that fixes the end $\xi_i$. Fix two reference points $o_1,o_2\in G_1,G_2$ respectively.

	For $i=1,2$, define the horocyclic function $h_i$ on $V(G_i)$ with respect to the end $\xi_i$ and reference point $o_i$ as follows:
	\[
	h_i(x_i)=\frac{\log \Delta(o_i,x_i)}{\log d_i}, x_i\in V(G_i),
	\]
	where $d_1=q,d_2=r$ and $\Delta(x,y)$ is the modular function for the subgroup $\Gamma_i$. (This definition differs by a negative sign as the one defined in some references like \cite{Bertacchi2001RW_on_DL_graphs}.)

 The Diestel-Leader graph $G=DL(q,r)$ consists of the couples $x_1x_2$ of $V(G_1)\times V(G_2)$ such that $h_1(x_1)+h_2(x_2)=0$, and  $x_1x_2$ is a neighbor of $y_1y_2$ if and only if $x_i$ is a neighbor of $y_i$ in $G_i$ for $i=1,2$. A schematic drawing $DL(2,2)$  can be found in Figure 2 on page 180 of \cite{Bartholdi_and_Woess2005spectral_computation_lamplighter_DL}. When $q\neq r$, the Diestel-Leader graph $G=DL(q,r)$ is a transitive nonunimodular graph. 
	
	The automorphism group $\mathrm{Aut}(G)$ of $G=DL(q,r)$ for $q\neq r$ can be described as 
	\[
	\mathrm{Aut}(G)=\{ \gamma_1\gamma_2\in \Gamma_1\times \Gamma_2\colon h_1(\gamma_1o_1)+h_2(\gamma_2o_2)=0  \};
	\]
	see \cite[Proposition 3.3]{Bertacchi2001RW_on_DL_graphs} for a proof. It is amenable since it is a closed subgroup of the amenable group $\Gamma_1\times \Gamma_2$. Hence Theorem \ref{thm: main theorem on a_n for amenable, nonuni} applies to Diestel-Leader graphs. Actually for Diestel-Leader graphs $G=DL(q,r),q\neq r$, it is known that  (\cite[Theorem 2]{Bartholdi_and_Woess2005spectral_computation_lamplighter_DL})
	\[
	u_{2n}\sim c_1\rho^{2n}\exp\big( -c_2n^{1/3} \big) n^{-5/6},
	\]
	where $\rho=\frac{2\sqrt{qr}}{q+r}$ is the spectral radius, and $c_1,c_2$ are explicit positive constants. 
\end{example}

%
%
%


Suppose $G_1,G_2$ are two transitive graphs with spectral radii $\rho_1,\rho_2$ and degrees $d_1,d_2$ respectively. It is well-known that the Cartesian product $G_1\times G_2$ has spectral radius $\rho=\rho(G_1\times G_2)=\frac{d_1\rho_1+d_2\rho_2}{d_1+d_2}$ (for instance see the proof of Proposition 18.1 in \cite{woess2000random}.) In fact the proof of Proposition 18.1 in \cite{woess2000random} also implies that if the return probabilities satisfy $u_n(G_i)\leq C_i\rho_i^n\cdot n^{\lambda_i}$ for some constants $C_i>0,i=1,2$, then the return probabilities on the Cartesian product satisfy $u_n(G_1\times G_2)\leq C\rho^n\cdot n^{\lambda_1+\lambda_2}$ for some constant $C>0$. 
\begin{example}[Cartesian product]\label{example: Cartesian product}
	Let $G_1$ be a connected, transitive graph. Let $G_2$ be a connected graph with a closed, amenable, nonunimodular, transitive subgroup $\Gamma$ of automorphisms. It is well known that the return probability on $G_1$ satisfies $u_n(G_1)\leq \rho(G_1)^n\cdot n^{\lambda_1}$ with $\lambda_1=0$ (for instance see (6.13) in \cite[Proposition 6.6]{LP2016}).  Theorem \ref{thm: main theorem on a_n for amenable, nonuni}  implies that the return probability on $G_2$ satisfies $u_n(G_2)\leq C\rho(G_2)^nn^{-3/2}$. Hence the above implication of the proof of Proposition 18.1 in \cite{woess2000random} yields that Conjecture \ref{conj: a_n is at most n^(-3/2)} also holds for the Cartesian product $G_1\times G_2$. 
\end{example}

\begin{example}[A free product]\label{example: free product}
	For  $G=C_\alpha*C_\beta\,(\beta,\alpha\geq 2,\max\{\alpha,\beta\}>2)$, the free product of two complete graphs of  $\alpha,\beta$ vertices respectively, one can show that $G$ has no closed, amenable and transitive subgroup. Actually if there is such a group $\Gamma$, then by Proposition \ref{prop: amenability equiv fix end} it must fix an end. But then it is easy to see that it can't be transitive. However this graph $G$ still has a closed, amenable, nonunimodular, \textbf{quasi-transitive} subgroup; see Example 4 on page 362 of \cite{Saloff-Coste_and_Woess1997transition_operators}. Hence the quasi-transitive case Theorem \ref{thm: main theorem on a_n for amenable, nonuni, quasi-transitive} applies. (Actually for such free products, $a_n\sim cn^{-3/2}$; see \cite{Woess1982local_limit_theorem_discrete_groups}.)
\end{example}

\section{Discussions on Conjecture \ref{conj: limit of the ratio of f_n over u_n}} \label{sec: discussion of Conj 1.6}

\subsection{A sufficient condition for Conjecture \ref{conj: limit of the ratio of f_n over u_n}}
 
Recall that for a connected transitive graph $G$ with spectral radius $\rho$, we denote by $u_n$ the 
$n$-step return probability for a simple random walk on $G$ and $a_n:=\frac{u_n}{\rho^n}$. 
\begin{proposition}\label{prop: a sufficient condition for Conjecture 4}
	Suppose $G$ is a locally finite, connected, transitive, transient graph. If for any $\varepsilon>0$, there exists $N=N(\varepsilon)>0$ such that for all  $n\geq 2N$ one has that
	\be\label{eq: sufficient condition for conj 4 in terms of un}
     \sum_{i=N}^{n-N}u_iu_{n-i}\leq \varepsilon u_n
	\ee
	or equivalently
	\be\label{eq: sufficient condition for conj 4 in terms of an}
	 \sum_{i=N}^{n-N}a_ia_{n-i}\leq \varepsilon a_n,
	\ee
	then Conjecture \ref{conj: limit of the ratio of f_n over u_n} holds for $G$. 
\end{proposition}
The inequality \eqref{eq: sufficient condition for conj 4 in terms of un} can be interpreted as conditioned on returning to the starting point at time $n$, the expected number of returns of the simple random walk to the starting point between time $N$ and $n-N$ is at most $\varepsilon$. Since $a_n=\frac{u_n}{\rho^n}$, the equivalence between \eqref{eq: sufficient condition for conj 4 in terms of un} and \eqref{eq: sufficient condition for conj 4 in terms of an} is obvious.

Before proving Proposition \ref{prop: a sufficient condition for Conjecture 4}, we first give some examples for Conjecture \ref{conj: limit of the ratio of f_n over u_n} using this proposition.

\subsection{Examples for Conjecture \ref{conj: limit of the ratio of f_n over u_n}}

\begin{lemma}\label{lem: some special cases for condition 6.1}
	 The condition \eqref{eq: sufficient condition for conj 4 in terms of un} holds if $(u_n)_N\geq0$ has one of the following asymptotic behavior:
	\begin{itemize}
		\item[(i)] $u_{2n}\asymp \rho^{2n}\cdot n^{-\alpha}$ for some constants $\alpha>1$ and $\rho\in(0,1]$,
		
		\item[(ii)]  $u_{2n}\asymp \rho^{2n}\cdot n^{-\alpha}\cdot e^{-cn^\beta}$ for some constants $\rho\in(0,1]$, $\alpha$ real, $c>0$ and $0<\beta<1$; 
		
		\item[(iii)] $u_{2n}\asymp \rho^{2n}\cdot e^{-n/(\log n)}$ for some constant $\rho\in(0,1]$. 
	\end{itemize} 
\end{lemma}

Lemma \ref{lem: some special cases for condition 6.1} is inspired  by Remark 1 in \cite{Chover_etal_1973}. If all odd terms $u_{2k+1}=0$, then one can verify condition \eqref{eq: sufficient condition for conj 4 in terms of un} easily in each of the three cases. If some odd terms $u_{2k+1}>0$, then by Lemma \ref{lem: ratio of u_n+1 over u_n tends to rho} the full sequence will satisfy the same asymptotic behavior instead of merely the even terms and hence condition \eqref{eq: sufficient condition for conj 4 in terms of un} can be verified similar to the  case of all odd terms being zero. We thus omit the details of the verification of Lemma \ref{lem: some special cases for condition 6.1}.

The reason for making Conjecture \ref{conj: limit of the ratio of f_n over u_n} is that there are a lot of examples support the conjecture. 

\begin{example}[graphs with polynomial growth rate]\label{example: polynomial growth graphs for Conj 4}
	If $G$ is a transient, transitive graph with polynomial growth rate, then  as discussed in Section \ref{sec: intro}, there exists an integer $k\geq3$ such that the volume of a ball with radius $n$ in $G$ has order $n^k$. Also the return probability satisfies $u_{2n}\asymp n^{-\frac{k}{2}}$; see Corollary 14.5, Theorem 14.12 and 14.19 in \cite{woess2000random}. Hence by Lemma \ref{lem: some special cases for condition 6.1} and Proposition \ref{prop: a sufficient condition for Conjecture 4}  such a graph $G$ satisfies Conjecture \ref{conj: limit of the ratio of f_n over u_n}. This was already noticed in \cite{Doney_and_Korshunov2011}. 
\end{example}

Conjecture \ref{conj: limit of the ratio of f_n over u_n} is open for general amenable Cayley graphs. For example we don't even know whether it holds for all Cayley graphs of certain lamplight groups; see the discussion after Example \ref{example: lamplighter graph for Conj 4}.

\begin{example}[hyperbolic graphs]\label{example: hyperbolic graphs for Conj 4}
	If $G$ is a hyperbolic graph, then one has that $a_{2n}\asymp n^{-3/2}$ \cite{Gouezel2014}. Hence \eqref{eq: sufficient condition for conj 4 in terms of an} is satisfied and then Conjecture \ref{conj: limit of the ratio of f_n over u_n} holds. This was already noticed by Gou\"{e}zel in \cite[Proposition 4.1]{Gouezel2014}. 
\end{example}

\begin{example}[free products]\label{example: free products for Conj 4}
	There are quite a lot  Cayley graphs of free products of groups  for which one knows well about the  asymptotic behavior of the return probabilities. We just mention a few of them here.
	\begin{itemize}
		\item[(i)] For the free products of two complete graphs as in Example \ref{example: free product}, one has that $a_{2n}\asymp n^{-3/2}$ by \cite{Woess1982local_limit_theorem_discrete_groups}.
		
		\item[(ii)] It was known that \cite{Cartwright_and_Soardi1986rw_free_products,Woess1986rw_free_produts} that the $n$-step return probabilities behaves like $u_{2n}\sim c\rho^{2n}n^{-3/2}$ under quite general conditions for random walks on a free product of discrete groups. For readers' convenience, quite a few of such conditions can be found in Corollary 6.12 of \cite{Woess1994survey}.

		\item[(iii)] For the free products  $\mathbb{Z}^d*\mathbb{Z}^d$ (natural generators, i.e., integer vectors   with Euclidean length one), one has that 
		\[
		a_{2n}\asymp 
		\left\{
		\begin{array}{cc}
			n^{-3/2} & \textnormal{ if }d\in\{1,2,3,4\}\\
			n^{-d/2} & \textnormal{ if }d\geq 5.
		\end{array}
		\right.
		\]
		This was due to Cartwright \cite{Cartwright1988RW_free_product}. Actually a general result holds for $\mathbb{Z}^d*\ldots*\mathbb{Z}^d$ ($s\geq 2$ times); see \cite{Cartwright1988RW_free_product} or \cite[Theorem 6.13]{Woess1994survey}. 
		
	\end{itemize}

	Given the explicit asymptotic behavior of return probabilities, it is easy to verify condition \eqref{eq: sufficient condition for conj 4 in terms of un} holds for all these examples and hence Conjecture \ref{conj: limit of the ratio of f_n over u_n} holds for them.
\end{example} 

 It seems to be new that the Examples \ref{example: lamplighter graph for Conj 4} and \ref{example: nonunimodular example for Conj 4} below satisfy Conjecture \ref{conj: limit of the ratio of f_n over u_n}. (Proposition 4.1 of \cite{Gouezel2014} also applies to graphs listed in Example \ref{example: free products for Conj 4}.)

\begin{example}[some Cayley graphs of lamplighter groups] \label{example: lamplighter graph for Conj 4}
	Consider a lamplighter group $H\wr\mathbb{Z}$, where $H$ is a finite group. 
	Revelle \cite[Theorem 1]{Revelle2003heat_kernel_lamplighter} showed that the return probability of simple random walk on the Cayley graph $G$ of the lamplighter group $H\wr \mathbb{Z}$ with a suitable chosen generating set satisfies 
	\[
	u_{2n}\sim c_2n^{1/6}\exp\big[-c_1n^{1/3}\big].
	\]
	Hence by Lemma \ref{lem: some special cases for condition 6.1} and Proposition \ref{prop: a sufficient condition for Conjecture 4} such a graph $G$ also satisfies Conjecture \ref{conj: limit of the ratio of f_n over u_n}. 
\end{example}

Unfortunately we don't even know whether Conjecture \ref{conj: limit of the ratio of f_n over u_n} hold for all Cayley graphs of such lamplighter group $H\wr \mathbb{Z}$. Theorem 1.1 of \cite{Pittet_Saloff_Coste2000stablity} says that if $\Gamma$ is a finitely generated group and $G_1,G_2$ are two Cayley graphs generated by symmetric finite generating sets of $\Gamma$, then the return probabilities on $G_1$ and $G_2$ satisfy
\[
u_n(G_1)\simeq u_n(G_2)
\]
in the sense that there exists a constant $C\geq 1$ so that 
\[
u_n(G_1)\leq C\cdot u_{n/C}(G_2) \textnormal{ and }u_n(G_2)\leq C\cdot u_{n/C}(G_1). 
\]
Applying this to Revelle's lamplighter group examples, one has that for any Cayley graph of $H\wr \mathbb{Z}$ the return probabilities satisfy
\[
c_4n^{1/6}\exp\big[-C_3n^{1/3}\big]\leq u_{2n}\leq C_3n^{1/6}\exp\big[-c_4n^{1/3}\big]
\]
for some constants $C_3,c_4>0$. However we are not able to verify \eqref{eq: sufficient condition for conj 4 in terms of un} with only this inequality. 

\begin{example}[some nonunimodular graphs]\label{example: nonunimodular example for Conj 4}
	As noted in Example \ref{example: grandparent graph}, for grandparent graph one has that $a_n\asymp n^{-3/2}$. Hence \eqref{eq: sufficient condition for conj 4 in terms of an} is satisfied and then Conjecture \ref{conj: limit of the ratio of f_n over u_n} holds by Proposition \ref{prop: a sufficient condition for Conjecture 4}.
	
	As noted in Example \ref{example: DL graph}, the explicit asymptotic behavior of return probabilities is known for Diestel--Leader graphs \cite[Theorem 2]{Bartholdi_and_Woess2005spectral_computation_lamplighter_DL}. By Lemma \ref{lem: some special cases for condition 6.1} and Proposition \ref{prop: a sufficient condition for Conjecture 4} one has that Conjecture \ref{conj: limit of the ratio of f_n over u_n} hold for all Diestel--Leader graphs. 
\end{example}

 In light of all these examples it is likely to be true  that the condition \eqref{eq: sufficient condition for conj 4 in terms of un} holds for all transient, transitive graphs.

\subsection{Proof of Proposition \ref{prop: a sufficient condition for Conjecture 4}}
A key ingredient for Proposition \ref{prop: a sufficient condition for Conjecture 4} is the following  theorem from \cite{Chover_etal_1973}. 
\begin{theorem}[Theorem 1 of \cite{Chover_etal_1973}]\label{thm: CNW1973}
	Let $\mu=\{\mu_n\}$ be a probability measure on nonnegative integers, where $\mu_n=\mu(n)$ is the mass of $n$. Let $r\geq1$ be the radius of the \notion{generating function} 
	\[
	\widehat{\mu}(z)=\sum_{n=0}^{\infty}\mu_nz^n.
	\]
	Assume that
	\begin{enumerate}
		\item[(i)] 
		\[
		\lim_{n\to\infty} \frac{\mu_n^{*2}}{\mu_n}:=\lim_{n\to\infty} \frac{\sum_{i=0}^{n}\mu_i\mu_{n-i}}{\mu_n}=C \,\,\textnormal{ exists }\,(<\infty);
		\] 
		
		\item[(ii)] 
		\[
		\lim_{n\to\infty}\frac{\mu_{n+1}}{\mu_n}=\frac{1}{r}\,\,(>0);
		\]
		\item[(iii)] $\widehat{\mu}$ converges at its radius of convergence:
		\[
		\widehat{\mu}(r)=D<\infty;
		\]
		
		\item[(iv)] 
		$\phi(w)$ is a function analytic in a region containing the range of $\widehat{\mu}(z)$ for $|z|\leq r$.
	\end{enumerate}
	
	Then there exists a measure $\phi(\mu)=\{\phi(\mu)_n,n\geq0\}$ on nonnegative integers with its generating function $\widehat{\phi(\mu)}(z):= \sum_{n=0}^{\infty}\phi(\mu)_nz^n$ satisfies
	\[
	\widehat{\phi(\mu)}(z)=\phi(\widehat{\mu}(z)), \textnormal{ for }|z|\leq r,
	\]
	and for which
	\be\label{eq: ratio of phi(mu)_n over mu_n}
	\lim_{n\to\infty}\frac{\phi(\mu)_n}{\mu_n}=\phi'(D).
	\ee
	
	
	Also we must have $C=2D$ in assumption (i).
\end{theorem}

%
%
%

 The following lemma is a special case $(x=y)$ of \cite[Theorem 5.2(b)]{Woess1994survey}.
 \begin{lemma}\label{lem: ratio of u_n+1 over u_n tends to rho}
 	Suppose $G$ is a locally finite, connected transitive graph with spectral radius $\rho$ and period $\mathsf{d}:=\gcd\{n\geq1,u_n>0 \}\in\{1,2\}$. Then 
 	\[
 	\lim_{n\to\infty,\mathsf{d}|n}\frac{u_{n+\mathsf{d}}}{u_n}=\rho^{\mathsf{d}}. 
 	\]
 \end{lemma}

 Recall that $U(z)=\sum_{n=0}^{\infty}u_nz^n$ and $F(z)=\sum_{n=0}^{\infty}f_nz^n$ are the generating functions associated with the return probabilities $(u_n)_{n\geq0}$ and first return probabilities $(f_n)_{n\geq0}$ respectively.  
  For recurrent, transitive graphs, the spectral radius $\rho$ satisfies $\rho=1$ and $U(1)=\infty, F(1)=1$. For transient, transitive graphs one has  the following simple result. 
 \begin{lemma}\label{lem: function U and F in the radius of convergence}
 	Suppose $G$ is a transient, transitive graph with spectral radius $\rho$.
 	Then 
 	\begin{enumerate}
 		\item[(a)] $U(\rho^{-1})<\infty$ and $F(\rho^{-1})<1$ and 
 		\item[(b)] for all complex number $z$ with $|z|\leq \rho^{-1}$, one has that 
 		\be\label{eq: relation of U(z) and F(z)} 
 		U(z)=\frac{1}{1-F(z)}. 
 		\ee
 	\end{enumerate}
 \end{lemma}
The inequality $U(\rho^{-1})<\infty$ in Part (a) of Lemma \ref{lem: function U and F in the radius of convergence} is just the fact $\sum_{n=0}^{\infty}a_n<\infty$ which we already mentioned in Section \ref{sec: intro} (Theorem 7.8 of \cite{woess2000random}). 
Part (b) of Lemma \ref{lem: function U and F in the radius of convergence} is basically contained in Lemma 1.13 of \cite{woess2000random} and then one can deduce $F(\rho^{-1})<1$ using $U(\rho^{-1})<\infty$ and \eqref{eq: relation of U(z) and F(z)}. The proof is thus omitted. 

\begin{proof}[Proof of Proposition \ref{prop: a sufficient condition for Conjecture 4}]
	Recall that $\mathsf{d}:=\gcd\{ n\geq 1\colon u_n> 0\}\in\{1,2\}$ denotes the period of simple random walk. We only deal with the case of $\mathsf{d}=1$; the case $\mathsf{d}=2$ is similar. 
	
	We shall use Theorem \ref{thm: CNW1973}. 
	In light of the relation \eqref{eq: relation of U(z) and F(z)} in Lemma \ref{lem: function U and F in the radius of convergence} it is natural to take the function $\phi:w\mapsto1-\frac{1}{U(1)w}$ and probability measure  $\mu=\{\mu_n,n\geq0\}$ given by $\mu_n=\frac{u_n}{U(1)},n\geq0$.
	Then 
	\[
	\widehat{\mu}(z)=\sum_{n=0}^{\infty}\mu_nz^n=\frac{U(z)}{U(1)} \textnormal{ has radius of convergence }r=\rho^{-1}.
	\]

	The assumption (i) in Theorem \ref{thm: CNW1973} now becomes 
	\be\label{eq: 2.2 to check for period =1 case}
	\lim_{n\to\infty} \frac{\mu_n^{*2}}{\mu_n}=\lim_{n\to\infty} \frac{\sum_{j=0}^{n}u_{j}u_{n-j}}{U(1)u_{n}}=C.
	\ee

	Assumption (ii) now becomes (and is verified by Lemma \ref{lem: ratio of u_n+1 over u_n tends to rho}):
	\be\label{eq: 4.5}
	\lim_{n\to\infty} \frac{\mu_{n+1}}{\mu_n}=\lim_{n\to\infty}\frac{u_{n+1}}{u_{n}}=\frac{1}{r}=\rho.
	\ee
	
	Assumption (iii) is also easy to verify in our set up:
	\be\label{eq: 4.6}
	\widehat{\mu}(r)=\frac{U(\rho^{-1})}{U(1)}=D<\infty.
	\ee
	
	As for assumption (iv), by Lemma \ref{lem: function U and F in the radius of convergence}
$
U(z)=\frac{1}{1-F(z)}
$ holds for all $|z|\leq \frac{1}{\rho}$. 
In particular $|U(z)|\geq \frac{1}{1-|F(\rho^{-1})|}>0$ for $|z|\leq \rho^{-1}$. Hence the function $\phi:w\mapsto 1-\frac{1}{U(1)w}$ is analytic in a region containing the range of  $\widehat{\mu}(z)=\frac{U(z)}{U(1)}$ for $|z|\leq r=\rho^{-1}$.
	
	The choice of $\phi$ yields that
	\[
	\widehat{\phi(\mu)}(z)=\phi(\widehat{\mu}(z))=1-\frac{1}{U(1)\widehat{\mu}(z)}=1-\frac{1}{U(z)}=F(z)=\sum_{n=1}^{\infty}f_{n}z^n, \textnormal{ for }|z|\leq r
	\]
	and 
	\[
	\phi'(D)=\frac{1}{U(1)D^2}\stackrel{\eqref{eq: 4.6}}{=}\frac{U(1)}{U(\rho^{-1})^2}.
	\]
	It is easy to see that if \eqref{eq: sufficient condition for conj 4 in terms of un} holds, then by \eqref{eq: 4.5} one has that \eqref{eq: 2.2 to check for period =1 case} holds for $C=2D=2\frac{U(\rho^{-1})}{U(1)}$.
	
	Hence if \eqref{eq: sufficient condition for conj 4 in terms of un} holds for a graph $G$, then all the assumptions of Theorem \ref{thm: CNW1973} hold. Thus one has that 
	\[
	\lim_{n\to\infty}\frac{\phi(\mu)_n}{\mu_n}\stackrel{\eqref{eq: ratio of phi(mu)_n over mu_n}}{=}\phi'(D)=\frac{U(1)}{U(\rho^{-1})^2}.
	\]
Since $\phi(\mu)_n=f_n$ and $\mu_n=\frac{u_n}{U(1)}$ one has that  Conjecture \ref{conj: limit of the ratio of f_n over u_n} holds for $G$:
	\[
	\lim_{n\to\infty} \frac{f_{n}}{u_{n}}=\frac{1}{U(\rho^{-1})^2}=[1-F(\rho^{-1})]^2. 
	\]

	If the period $\mathsf{d}=2$, it is easy to see that $u_{2n+1}=0$ for all $n$. Hence we just take the probability measure $\mu=\{\mu_n,n\geq0\}$ to be given by $\mu_n=\frac{u_{2n}}{U(1)},n\geq0$. In this case $r=\rho^{-2}$ and $\widehat{\mu}(z)=\frac{U(\sqrt{z})}{U(1)}$ for $|z|\leq \rho^{-2}$. The rest is similar to the case of $\mathsf{d}=1$ and omitted.
\end{proof}

\subsection{Final remark about Conjecture \ref{conj: limit of the ratio of f_n over u_n}}

Recall that condition \eqref{eq: sufficient condition for conj 4 in terms of un} roughly says that conditioned on returning to the starting point at time $n$, the \textbf{expectation} of returns of the simple random walk between time $N$ and $n-N$ is small for large $N$. Proposition \ref{prop: a sufficient condition for Conjecture 4} says that if  \eqref{eq: sufficient condition for conj 4 in terms of un} holds, then Conjecture \ref{conj: limit of the ratio of f_n over u_n} holds. We remark that on the other hand if Conjecture \ref{conj: limit of the ratio of f_n over u_n} holds, then  Conjecture \ref{conj: one big excursion} holds. Here Conjecture \ref{conj: one big excursion} roughly says that conditioned on returning to the starting point at time $n$, \textbf{with high probability} most of the returns of the simple random walk  occurred near time $0$ or $n$. 


Suppose $G=(V,E)$ is a locally finite, connected, transitive, transient graph with spectral radius $\rho$. Fix an arbitrary vertex $o\in V$.  Let $(X_n)_{n\geq0}$ be a simple random walk on $G$ starting from $o$. Write $f_n$ for the first return probability at time $n$ and $F(z)=\sum_{n=1}^{\infty}f_nz^n$ for the corresponding generating function. Let $\mathsf{d}$ denote the period of the simple random walk. We will consider the returning times to $o$ conditioned on $\{X_n=X_0=o\}$. Define the returning times $(s_i)_{i\geq0},(l_i)_{i\geq0}$ as follows (here $(l_i)_{i\geq0}$ records the returning times in the reverse order):
\begin{itemize}
	\item $s_0=l_0=0$ and ,
	\item for $i\geq 0$, 
	\[
	s_{i+1}=\min\big\{k\colon k>s_i,X_{k}=o  \big\},\,\,\,l_{i+1}=\min\big\{k\colon  k>l_i,X_{n-k}=o \big\}.
	\]
\end{itemize}
Let 
\[
\alpha=\alpha(n)=\max\big\{  k\geq0\colon s_k\leq \frac{n}{2}\big\},\beta=\beta(n)=\max\big\{ k\geq0\colon l_k\leq \frac{n}{2} \big\}. 
\]
Consider the random variable 
\[
\mathsf{V}_{n}=\big( (s_1,\ldots,s_\alpha,0,0,\ldots),(l_1,\ldots,l_\beta,0,0,\ldots) \big)
\]
which takes values in the space $\mathbb{N}^{\mathbb{N}}\times \mathbb{N}^{\mathbb{N}}$.
\begin{conj}\label{conj: one big excursion}
	For any transient, transitive graph, the distribution of $\mathsf{V}_{n} $ conditioned on the event $\{X_n=X_0=o\}$ converges as  $n\to\infty,\mathsf{d}|n $, to the distribution of the random variable 
	\[
	\big( (T_1,\ldots,T_L,0,0,\ldots),(\hat{T}_1,\ldots,\hat{T}_{\hat{L}},0,0,\ldots) \big)
	\]
	where $(T_j)_{j\geq1}$ are the partial sums of an i.i.d.\ sequence $(\xi_i)_{i\geq1}$ with distribution given by  $\mathbb{P}[\xi_i=k]=\frac{f_{k}\rho^{-k}}{F(\rho^{-1})}$ and $L$ is an independent random variable  with a geometric distribution  with parameter $1-F(\rho^{-1})$, and  $(\hat{T}_1,\ldots,\hat{T}_{\hat{L}},0,0,\ldots)$ is an independent copy of 
	$(T_1,\ldots,T_L,0,0,\ldots)$.  
	
\end{conj}

Conjecture \ref{conj: one big excursion} is inspired by   \cite[Proposition 2.2]{BJ1999brownian_bridge} which says that Conjecture \ref{conj: one big excursion} holds for regular trees. The sketch below is also a simple modification of the proof of \cite[Proposition 2.2]{BJ1999brownian_bridge}.

\begin{proof}[Sketch of the implication of Conjecture \ref{conj: limit of the ratio of f_n over u_n} $\Rightarrow$ Conjecture \ref{conj: one big excursion}]
We only deal with the case $\mathsf{d}=1$ here; the case of $\mathsf{d}=2$ can be treated similarly.	
If Conjecture \ref{conj: limit of the ratio of f_n over u_n} holds and $\mathsf{d}=1$, then
	\[
	\lim_{n\to \infty}\frac{f_n}{u_n}=\big(1-F(\rho^{-1})\big)^2\in(0,1).
	\]
	
When $m$ is fixed and $n\to\infty$, by the above limit and Lemma \ref{lem: ratio of u_n+1 over u_n tends to rho} one has that 
\begin{equation*}
\frac{f_{n-m}}{u_n}=\frac{f_{n-m}}{u_{n-m}}\cdot \frac{u_{n-m}}{u_n}\sim \big(1-F(\rho^{-1})\big)^2\cdot \rho^{-m}. 
\end{equation*}	
Therefore when $n$ is large, if $m=\sum_{i=1}^{a}k_i+\sum_{j=1}^{b}r_b$, then 
\begin{eqnarray}\label{eq: asymptotic behavior of Vn}
&&	\mathbb{P}\Big[ \alpha=a,s_i=\sum_{t=1}^{i}k_t,i\in\{1,\ldots,a\},\beta=b, l_j=\sum_{t=1}^{j}r_t,j\in\{1,\ldots,b\}\mid X_n=X_0=o \Big]\nonumber\\
&=&\bigg(\prod_{i=1}^{a}f_{k_i}\bigg)\cdot \frac{f_{n-m}}{u_n}\cdot\bigg( \prod_{j=1}^{b}f_{r_j}\bigg)\nonumber\\
&\sim& \bigg(\prod_{i=1}^{a}f_{k_i}\bigg)\cdot\bigg( \prod_{j=1}^{b}f_{r_j}\bigg)\cdot \big(1-F(\rho^{-1})\big)^2\cdot \rho^{-m}\nonumber\\
&=&\bigg(\prod_{i=1}^{a}f_{k_i}\rho^{-k_i}\bigg)\cdot\bigg( \prod_{j=1}^{b}f_{r_j}\rho^{-r_j}\bigg)\cdot \big(1-F(\rho^{-1})\big)^2,
\end{eqnarray}	
where we use the convention that $\prod_{i=1}^{0}=1$. 
Note that the last expression in \eqref{eq: asymptotic behavior of Vn} gives a probability measure since
\[
\sum_{a\geq0,k_i\geq1,b\geq0,r_j\geq1}\bigg(\prod_{i=1}^{a}f_{k_i}\rho^{-k_i}\bigg)\cdot\bigg( \prod_{j=1}^{b}f_{r_j}\rho^{-r_j}\bigg)\cdot \big(1-F(\rho^{-1})\big)^2=1. 
\]
From \eqref{eq: asymptotic behavior of Vn} it is easy to obtain the desired conclusion; for instance to see the distribution of $\alpha$ is tending to Geometric with parameter $1-F(\rho^{-1})$, it suffices to sum \eqref{eq: asymptotic behavior of Vn} over all possible $k_i,r_j,b$. 
\end{proof}

\section*{Acknowledgment}
The work is supported by ERC starting grant 676970 RANDGEOM. We thank Asaf Nachmias for informing Conjecture \ref{conj: a_n is at most n^(-3/2)} and many helpful discussions. We thank Russ Lyons for pointing out the proof of  Claim \ref{claim: same radii of conv for F(z) and U(z) for graphs} using Pringsheim's theorem. We thank two referees for their careful reading and helpful comments.

\appendix
\section{Proof of Lemma \ref{lem: prob of M_n=rt_0, Y_n=0 for quasi-transitive}}


\begin{proof}[Proof of Lemma \ref{lem: prob of M_n=rt_0, Y_n=0 for quasi-transitive}]
	Let $\Omega=\{(i,j,q)\colon  N_{i,j,q}\neq \emptyset \}$. Let $(\xi_n)_{n\geq0}$ be a Markov chain on $\Omega$ induced by the $p_h$-walk $(S_n)_{n\geq0}$. More precisely, the initial distribution  of $\xi_1$ is given by 
	\[
	\mathbb{P}[\xi_1=(i,j,q)]=\mathbb{P}[S_0=o_i,S_1\in\Gamma o_j,\Delta(S_0,S_1)=q]=\big(\sqrt{q}\#N_{i,j,q}\big)\cdot (v_iv_j)\cdot \frac{1}{\rho\sqrt{d_id_j}},
	\]
	and the transition probability is given by 
	\begin{eqnarray*}
		\mathbb{P}\big[ \xi_{k+1}=(i',j',q') \mid \xi_{k}=(i,j,q) \big]
		&=&\mathbf{1}_{\{ i'=j \}}\cdot \mathbb{P}[ S_{k+1}\in \Gamma o_{j'} \textnormal{ and }\Delta(S_k,S_{k+1})=q' \mid S_k\in \Gamma o_j]\nonumber\\
		&=&\mathbf{1}_{\{ i'=j \}}\cdot \big(\sqrt{q'}\#N_{j,j',q'}\big)\cdot \frac{v_{j'}}{v_j}\cdot \frac{1}{\rho\sqrt{d_jd_{j'}}},
	\end{eqnarray*}
	Obviously $(\xi_n)_{n\geq1}$ is a finite, irreducible Markov chain starting from the stationary probability measure. 
	
	Let $f:\Omega\to \mathbb{R}$ be a function defined by $f\big((i,j,q)\big)=\log q$. Write $Z_k=f(\xi_k)$ for $k\geq 1$. Then it is easy to see that $(Y_n)_{n\geq0}$ has the same law as  the partial sums of the sequence $(Z_k)_{k\geq1}$. So in the following we will assume that $(Y_n)_{n\geq0}$ are the partial sums: $Y_0=0,Y_n=\sum_{k=1}^{n}Z_k$ for $n\geq1$. 
	
	The first step is to prove the  ballot theorem in this setup.
	\begin{claim}\label{claim: ballot theorem for quasi-transitive}
		There is a constant $c>0$ such that for all $0\leq k\leq n$, 
		\be\label{eq: ballot theorem for quasi-transitive}
		\mathbb{P}\big[Y_j>0,j=1,\cdots,n-1,Y_n\in[kt_0,(k+1)t_0)\big]\leq c\frac{k\vee 1}{n^{3/2}}.
		\ee
	\end{claim}
	\begin{proof}[Proof of Claim \ref{claim: ballot theorem for quasi-transitive}]
		We follow the proof of Theorem 1 in \cite{addario2008ballot}.
		
		First by Theorem 1 in \cite{Bolthausen1980} there is a constant $c_1>0$  such that for all $n$,
		\be\label{eq: concentration bounds}
		\sup_{x\in\mathbb{R}}\mathbb{P}[x\leq Y_n\leq x+t_0]\leq \frac{c_1}{\sqrt{n}}.
		\ee

		 Secondly we show that the item (iii) of Lemma 3 in \cite{addario2008ballot} still holds in this setup, namely,  for $h\geq0$ and $T_h(Y):=\inf\{n\colon Y_n<-h\}$,
		\be\label{eq: A.3}
		\mathbb{P}[T_h(Y)\geq n]\leq c\frac{h\vee 1}{\sqrt{n}}. 
		\ee
		
		Fix an arbitrary $x=(i,j,q)\in\Omega$ and write $\mathbb{P}_x,\mathbb{E}_x$ for the law of the Markov chain $(\xi_n)_{n\geq1}$ and expectation conditioned on $\xi_1=x$. Also  let $R_k$ be the $k$-th return to $x$ of the Markov chain $(\xi_n)_{n\geq0}$, i.e., 
		$R_1=\inf\{k\geq 1 \colon \xi_k=x \}$ and $R_n=\inf\{k> R_{n-1} \colon \xi_k=x \}$ for $n\geq 2$. Let $U_i=\sum_{k=R_i}^{R_{i+1}-1}Z_k$ be the sum of the $i$-th excursion.  Since during each excursion the expected number of visits to the states $y\in \Omega$ is a stationary measure (see Theorem 6.5.2 of \cite{Durrett2010Probability_examples_4th_edition}), one has $\mathbb{E}_x[U_i]=0$ by \eqref{eq: mean zero when start with a stationary vertex}. Hence $(U_i)$ are i.i.d.\ r.v.'s with mean zero. Let $\Lambda_n=\max\{k\colon R_k\leq n\}$ be the number of returns to $x$ up to time $n$. By a large deviation principle, for $\beta=2\mathbb{E}[R_2-R_1]$ there is a constant $c_2>0$ such that
		\be
		\mathbb{P}_x[\Lambda_n \leq \frac{n}{\beta}]\leq \frac{\exp(-c_2n)}{c_2}.
		\ee
		For $h\geq0$,  let $T_h(U)=\inf\{n\colon \sum_{i=1}^{n}U_i<-h\}$. Therefore
		\begin{eqnarray}
			\mathbb{P}_x[T_h(Y)\geq n]&\leq&
			\mathbb{P}_x[\Lambda_n \leq \frac{n}{\beta}]+ \mathbb{P}_x[\Lambda_n>\frac{n}{\beta},T_h(Y)\geq n]\nonumber\\
			&\leq&
			\mathbb{P}_x[\Lambda_n \leq \frac{n}{\beta}]+ \mathbb{P}_x[T_h(U)\geq \frac{n}{\beta}]\nonumber\\
			&\leq& \frac{\exp(-c_2n)}{c_2}+\frac{c_3(h\vee 1)}{\sqrt{n/\beta}}\leq c_x\frac{h\vee 1}{\sqrt{n}}
		\end{eqnarray}
		where in the last step we use the item (iii) of Lemma 3 in \cite{addario2008ballot} for the i.i.d.\ sequence $(U_i)$. Taking $c=\max\{c_x\colon x\in\Omega\}$ one has \eqref{eq: A.3}. 
		
		Since $(\xi_n)_{n\geq 0}$ is an irreducible Markov chain with a finite state space $\Omega$, there exists a constant $\delta>0$ such that for any $x,y\in\Omega, n\geq 1$, if $\mathbb{P}[\xi_n=y\mid \xi_1=x]>0$, then  $\mathbb{P}[\xi_n=y\mid \xi_1=x]>\delta$. Hence for any $x,y\in\Omega$ such that $\mathbb{P}[\xi_n=y\mid \xi_1=x]>0$, by \eqref{eq: concentration bounds} one has that 
		\be\label{eq: concentration bounds when conditioned on endpoints}
		\sup_{t\in\mathbb{R}}\mathbb{P}[t\leq Y_n\leq t+t_0\mid \xi_1=x,\xi_n=y]\leq \frac{c_1}{\mathbb{P}[\xi_1=x]\delta\sqrt{n}}= \frac{c_4}{\sqrt{n}}.
		\ee
		Similarly  for any $x,y\in\Omega$ such that $\mathbb{P}[\xi_n=y\mid \xi_1=x]>0$,
		\be\label{eq: tail probability for first hitting time of a level h}
		\mathbb{P}_x[T_h(Y)\geq n\mid \xi_n=y]\leq c\frac{h\vee 1}{\sqrt{n}}.
		\ee
		
		Now fix a pair $x,y\in\Omega$ such that $\mathbb{P}[ \xi_{\lfloor \frac{n}{4}\rfloor}=x,\xi_{\lceil \frac{3n}{4}\rceil}=y]>0$. Consider the probability 
		\[
		L_{k,n}=L_{k,n}(x,y):=\mathbb{P}\big[Y_j>0,j=1,\cdots,n-1,Y_n\in[kt_0,(k+1)t_0), \xi_{\lfloor \frac{n}{4}\rfloor}=x,\xi_{\lceil \frac{3n}{4}\rceil}=y\big].
		\]
		Let $Y^r$ be the sequence given by  $Y^r_0=0$ and for $i$ with $0\leq i<n$, $Y^r_{i+1}=Y^r_i-Z_{n-i}$, i.e., partial sums of the sequence $(-Z_{n-i})_{i=0}^{n-1}$. For $h\geq0$, let $T_h^r(Y)$ be the minimum of $n$ and the first time $t$ that $Y^r_t\leq -h$. By considering the reversed chain of $(\xi_n)_{n\geq1}$ and $-f$, one has that \eqref{eq: tail probability for first hitting time of a level h} also holds for $T_h^r(Y)$, in particular, 
		\be\label{eq: tail probability for first hitting time of a level h, reversed chain}
		\mathbb{P}\big[T^r_{(k+1)t_0}(Y)>\lfloor \frac{n}{4}\rfloor \mid \xi_{\lceil \frac{3n}{4}\rceil}=y \big]\leq c\frac{(k+1)t_0\vee 1}{\sqrt{n}}.
		\ee

		In order that $Y_n\in[kt_0,(k+1)t_0) $ and $Y_i>0$ for all $0<i<n$, it is necessary that 
		\begin{enumerate}
			\item[(a)] $T_0(Y)>\lfloor \frac{n}{4}\rfloor$,
			
			\item[(b)]   $T^r_{(k+1)t_0}(Y)>\lfloor \frac{n}{4}\rfloor$, and 
			
			\item[(c)] $Y_n\in [kt_0,(k+1)t_0)$.
		\end{enumerate} 
		Writing $g_{k,n}(x,y)=\mathbb{P}\big[kt_0\leq Y_n<(k+1)t_0 \mid T_0(Y)>\lfloor \frac{n}{4}\rfloor, T^r_{(k+1)t_0}(Y)>\lfloor \frac{n}{4}\rfloor, \xi_{\lfloor \frac{n}{4}\rfloor}=x,\xi_{\lceil \frac{3n}{4}\rceil}=y\big]$, one has that
		\begin{eqnarray}\label{eq: 4.12}
			L_{k,n}&\leq&\mathbb{P}\big[T_0(Y)>\lfloor \frac{n}{4}\rfloor, T^r_{(k+1)t_0}(Y)>\lfloor \frac{n}{4}\rfloor, \xi_{\lfloor \frac{n}{4}\rfloor}=x,\xi_{\lceil \frac{3n}{4}\rceil}=y \big]\cdot g_{k,n}(x,y)\nonumber\\
			&=& \mathbb{P}[\xi_{\lfloor \frac{n}{4}\rfloor}=x,\xi_{\lceil \frac{3n}{4}\rceil}=y]\cdot 
			\mathbb{P}\big[T_0(Y)>\lfloor \frac{n}{4}\rfloor \mid \xi_{\lfloor \frac{n}{4}\rfloor}=x\big]\cdot \nonumber\\
			&& \mathbb{P}\big[T^r_{(k+1)t_0}(Y)>\lfloor \frac{n}{4}\rfloor \mid \xi_{\lceil \frac{3n}{4}\rceil}=y \big]\cdot g_{k,n}(x,y)\nonumber\\
			&\stackrel{\eqref{eq: tail probability for first hitting time of a level h},\eqref{eq: tail probability for first hitting time of a level h, reversed chain}}{\leq}& c^2\frac{(k+1)t_0\vee 1}{n}\cdot \mathbb{P}[\xi_{\lfloor \frac{n}{4}\rfloor}=x,\xi_{\lceil \frac{3n}{4}\rceil}=y]\cdot g_{k,n}(x,y)
		\end{eqnarray}
	where in the second step we use Markov property for $(\xi_n)_{n\geq0}$. By Markov property and \eqref{eq: concentration bounds when conditioned on endpoints} (applied to $Y_{\lceil \frac{3n}{4}\rceil}-Y_{\lfloor \frac{n}{4}\rfloor}$ conditioned on $\xi_{\lfloor \frac{n}{4}\rfloor},\xi_{\lceil \frac{3n}{4}\rceil},Y_{\lfloor \frac{n}{4}\rfloor}$ and $Y^r_{\lfloor \frac{n}{4}\rfloor}$)  one has that $g_{k,n}(x,y)\leq \frac{c_4}{\sqrt{n/2}}$. Therefore summing \eqref{eq: 4.12}  over all possible pairs $(x,y)\in\Omega\times \Omega$ such that $\mathbb{P}[ \xi_{\lfloor \frac{n}{4}\rfloor}=x,\xi_{\lceil \frac{3n}{4}\rceil}=y]>0$, one has \eqref{eq: ballot theorem for quasi-transitive}.  
	\end{proof}
	
	The second step is show the following analogue of Lemma \ref{lem: decay of hitting probability of level r}:
	\begin{claim}\label{claim: decay of tau_r=k}
		Let $\tau_r:=\inf\{ i\geq0\colon Y_i\geq rt_0 \}.$ One has that for $r\geq 1$,
		\be\label{eq: an upper bound on tau_r=k in the quasi-transitive}
		\mathbb{P}[\tau_r=k]\preceq \frac{r}{k^{3/2}}.
		\ee	
	\end{claim}
	\begin{proof}[Proof of Claim \ref{claim: decay of tau_r=k}]
		Consider the reversed chain $(\widetilde{\xi}_n)_{n\geq0}$ of $(\xi_n)_{n\geq0}$ started from the stationary distribution. Let $\widetilde{Z}_k=f(\widetilde{\xi}_n)$. 
		The vector 	 $(Z_n,\cdots, Z_1)$ has the same distribution as $(\widetilde{Z}_1,\widetilde{Z}_{2},\cdots,\widetilde{Z}_n)$. Let $\widetilde{Y}_n$ be the partial sums of $(\widetilde{Z}_k)_{k\geq 1}$. 
		
		The rest is the same as the proof of Lemma \ref{lem: decay of hitting probability of level r} just by replacing the ballot theorem by Claim \ref{claim: ballot theorem for quasi-transitive} for $\widetilde{Y}$ instead. 
	\end{proof}
	
	Now we are ready to show \eqref{eq: return to level 0 with maximum in level r for quasi-transitive}.

		Similar to  \eqref{eq: ballot thm adapted form}, for $k\in[r,n-r]$ using Markov property and Claim \ref{claim: ballot theorem for quasi-transitive} one has that 
		\be
		\mathbb{P}\big[Y_j-Y_k<t_0,j=k+1,\cdots, n,  Y_n-Y_k \in(-(r+1)t_0,-rt_0] \mid Y_k,\xi_k\big]  \leq c\frac{r+1}{(n-k)^{3/2}}.
		\ee
		Taking expectation one has that
		\be
		\mathbb{P}\big[M_n\in[rt_0,(r+1)t_0), Y_n=0,\tau_r=k\big]
		\leq \mathbb{P}[\tau_r=k]\cdot c\frac{r+1}{(n-k)^{3/2}}.
		\ee
		Hence similar to the deduction of \eqref{eq: 2.14}, we have \eqref{eq: return to level 0 with maximum in level r for quasi-transitive}.
\end{proof}

\bibliography{re_prob_nonuni_ref}
\bibliographystyle{plain}

\end{document}